\documentclass[10pt]{article}

\usepackage{graphicx}
\usepackage{amssymb,amsfonts,amsmath}
\usepackage[english]{babel}
\usepackage{pstricks}
\usepackage{longtable}
\usepackage{subfigure}
\usepackage{array}
  \usepackage{esint}    

\newtheorem{Theorem}{Theorem}
\newtheorem{Proposition}[Theorem]{Proposition}

\newtheorem{Lemma}[Theorem]{Lemma}

\begin{document}

\noindent
\section*{Existence and Uniqueness of Solutions for Coupled Hybrid Systems of Differential Equations}

\begin{center}
M. Menci and M. Papi \footnote{Universit\'a Campus Bio-Medico di Roma, Department of Engineering, Via \'Alvaro del Portillo 21, 00128, Roma, Italy.}
\end{center}

\begin{abstract}
In this paper we propose local and global existence results for the solution of systems characterized by the coupling of ODEs and PDEs. The coexistence of distinct mathematical formalisms represents the main feature of \textit{hybrid} approaches, in which the dynamics of interacting agents are driven by second-order ODEs, while reaction-diffusion equations are used to model the time evolution of a signal influencing them. We first present an existence result of the solution, locally in time.
In particular, we generalize the framework of recent works presented in the literature, concerning collective motions of cells due to mechanical forces and chemotaxis, taking into account a uniformly parabolic operator with space-and-time dependent coefficients, and a more general structure for the equations of motion. 
Then, the previous result is extended in order to obtain a global solution. \\

\end{abstract}

\noindent
\textbf{Keywords}
Hybrid Models, Parabolic Differential Equation,  Collective Motions, Chemotaxis
%


\section{Introduction}
\label{intro}

In this paper we study a particular type of systems of differential equations, arising from mathematical models that simultaneously combine discrete and continuum approaches, known as \textit{hybrid models}.
In recent years, these kinds of models have been mainly adopted to describe phenomena concerning living systems, such as cell aggregates or crowds (\cite{anderson}, \cite{cardio}, \cite{zebrafish}, \cite{colombi}, \cite{cristiani}). These are regarded as collections of agents presenting two fundamental aspects: a proper behavior and the ability to sense and actively interact with other individuals and the surrounding environment. 
The majority of mathematical models in the literature treats agents aggregates either as a continuum or as a discrete set of individuals.


Discrete models operate at the scale of individuals (\cite{berec}, \cite{cucker}, \cite{dormann}, \cite{wootton}). For example, in the biomedical field, each cell is treated as a unit of finite volume, which is able to move, divide and die individually according to biological observations. 
Agents have been modelled with simple points, spheres and ellipsoids, both with fixed volume and size, or more complex evolving deformable structures (\cite{dallon},\cite{palsson}).
A discrete approach allows to easily model mechanical interactions with other cells and with the surrounding matrix, and to incorporate details concerning individual cells (e.g. size, metabolic state). 
The possibility to model agents in fine details leads to the drawback of a large computational cost, which rapidly increases with the number of agents considered. 
Thus, discrete models are suitable for a microscopic description of phenomena when the number of agents is relatively limited. 

Since a same problem can be modelled at different scales depending on the aspects one is interested in, the choice of approach to adopt is often not unique.
Sometimes the distinction among different approaches is not extremely defined, other times it is quite obvious. In fact, to model regions in which mechanical and rheological properties are of primary interest, it is completely unnecessary to focus of the cell scale.
Aiming at a global description of the agents interactions, continuum models better fulfil the requirements (\cite{chaplain},\cite{colim}).
From a macroscopic point of view, the units' aggregate is described through its spatial mass density, and continuum equations are used to model cell-cell and cell-matrix interactions.
On the one hand, continuous models are easier to be computationally analyzed, and do not present limitations in the number of agents involved. On the other hand, they suffer from the fact that generally the averaging over space realized in continuum formalisms cannot fully account for the diversity of cellular and sub-cellular dynamical features.\\
Advantages and disadvantages of the two categories seem to be complementary.
The emergent hybrid approaches, in which some of the variables of the model are continuous and other are discrete, gain the advantages offered by both, providing a link between macro and microscale descriptions. 
In the literature, there is no a commonly agreed definition for the class of hybrid models.
In \cite{anderson}, authors investigate the effects of individual-based cell interactions in the different stages of tumor growth, presenting a model defined by a system of coupled non-linear partial differential equations.  

As already mentioned, modelling certain cell processes with a pure continuum approach would be challenging. The authors classify their model as \textit{hybrid}, since they discretize the reaction-diffusion model using finite-difference methods in order to focus on the individual cell level. In particular, they consider a random-walk model, assuming that the coefficients arising from the discretization process correspond to the probabilities for cells to move or remaining at their current location.
Authors refer to this kind of procedure as \textit{hybrid discrete-continuum}. 
This technique has been used in other works related to cancer growth \cite{anderson2}, and also in the context of angiogenesis \cite{angio} and retinal vasculature development \cite{retina}.

In our work we consider a hybrid differential system inspired by the mathematical models presented in \cite{cardio} and \cite{zebrafish}. 
In those papers, cells are treated as a set of localized agents whereas chemical concentrations are described through the spatial distribution of their concentrations.
%
The main difference is that in \cite{cardio}, to describe proliferation and differentiation occurring in the examined phenomenon, stochasticity is introduced in the model.
In the above mentioned works, authors focus on the construction of the models and provide some numerical simulations able to reproduce the related biological phenomena. 
To this end, only a numerical approximation of the solution of the system occurring in the mathematical model has been proposed.
However, while applications of hybrid models are increasingly frequent, the literature concerning wellposedness (existence and uniqueness) of the solution of the resulting differential systems is still lacking.
  
This is why our results can be considered as a first attempt to formalize, validate and justify the emergent class of hybrid models, in particular the ones in the form considered in \cite{cardio} and \cite{zebrafish}, from a theoretical point of view.

From a mathematical point of view, these models share some common features. The dynamic of the cells is driven by second-order ODEs, whereas the evolution in time of the chemical signals is described by reaction-diffusion equations.
Moreover, the dynamics of agents are not only influenced by mechanical interactions, but also by \textit{chemotaxis}, which leads the agents from regions with low concentration to the ones at an higher level.
In our model we consider a group of $n$ interacting agents. The dynamic of each agent $i=1,...,n$ is expressed by the following second-order differential equation:
\begin{equation}\label{intro1}
\mathbf{\ddot{x}}_i(t) = \displaystyle \mathbf{F}_{i}\left(t, \mathbf{X}(t),  \mathbf{\dot{X}}(t), \nabla f \left( \mathbf{x}_i(t),t;  \mathbf{X} \right)\right),
\end{equation}
where $\mathbf{X}(t)=\left[\mathbf{x}_1(t),...,\mathbf{x}_n(t)\right], \mathbf{\dot{X}}(t)=\left[\mathbf{\dot{x}}_1(t),...,\mathbf{\dot{x}}_n(t)\right]$  are the collections of position and velocity of each agents, at each time $t $, and $f$ models a signal influencing the dynamics (e.g. the concentration of a chemoattractant as in \cite{zebrafish}, \cite{cardio}). The notation $f=f(x,t;\mathbf{X} )$ aims to emphasize the dependence on the whole trajectory of all agents till time $t$, that will be clear in the next section (see equation ($\ref{fgamma}$)) . \\
In the proposed model, $f$ is the solution of the Cauchy problem 
\begin{equation}\label{intro2}
\left\{
\begin{array}{rl}
Lf(x,t;  \mathbf{X} )&=g(x,\mathbf{X}(t)), \ \ \ \ \ \ \ \ (x,t) \in \mathbb{R}^{N} \times (0,T)\\\
f(x,0)&= \varphi(x) \ \ \ \ \ \ \ \ \ \ \ \ \ \ \ \ \ x \in \mathbb{R}^{N},
\end{array}
\right.
\end{equation}
where $L$ is the following differential operator:
\begin{equation}\label{intro3}
L = \displaystyle \sum_{i,j=1}^{N} a_{i,j}(x,t) \partial_{x_i,x_j}^{2}+ 
	 \sum_{i=1}^{N} b_{i}(x,t) \partial_{x_i}+c(x,t)- \partial_{t}.
\end{equation}

Moreover we investigate a different variation of (\ref{intro1}), in which $\nabla f (\mathbf{x}_i (t), t; \mathbf{X})$  is replaced by the average over a ball centered in $\mathbf{x}_i (t)$ and having radius $\delta>0$. The introduction of an average gradient in \cite{cardio} and \cite{zebrafish} aims at modeling the fact that, in biological system, a cell feels the presence of chemical signals not only in its center, but also in the region surrounding it.\\
Although the operator $L$ is linear, we remark that, to the best of our knowledge, in the literature the mathematical foundations of such a problem have not been investigated. Therefore we follow a step-by-step approach to the problem, leaving the non linear case for a future work. 

The paper is organized as follow: in Section 2 we present the system of differential equations (\ref{intro1})-(\ref{intro3}) that will be considered throughout the whole paper. In particular, we will state the main assumptions under which the results will be given. In Section 3 we prove a local existence and uniqueness theorem, using a fixed point argumentation, even considering a slight modification of ($\ref{intro2}$). In Section 4, using a principle of continuation of the solution, the previous results are extended to global ones. 
Conclusive remarks and future works directions are collected in Section 6.
A short appendix completes the paper, followed by a glossary section, in which we describe the notation used throughout the paper.

\ \ \ \


\section{Problem statement}
\label{sec:1}

We consider the following system of differential equations
\begin{equation}\label{dinamica}
\left\{
\begin{array}{rl}
\mathbf{\dot{x}}_i(t)&=\mathbf{v}_{i}(t), \ \ \ \ \mbox{$0<t<T$},\\\
\mathbf{\dot{v}}_i(t)&= \displaystyle \mathbf{F}_{i}\left( t, \mathbf{X}(t), \mathbf{V}(t), \nabla f \left( \mathbf{x}_i(t),t;  \mathbf{X} \right)\right), \ \ \ \  i=1,...,n,
\end{array}
\right.
\end{equation}
with initial data $
\mathbf{x}_{i}(0)=\mathbf{x_{0}}_{i}, \ \ \ \ \mathbf{v}_{i}(0)=\mathbf{v_{0}}_{i} \ \ \ \in \mathbb{R}^N, \ \ \ \ \forall i=1,...,n.
$
Here
\begin{equation*}\label{hp_F}
\mathbf{F}_i:  \left[0,T\right] \times \mathbb{R}^{N\times n} \times \mathbb{R}^{N\times n} \times  \mathbb{R}^{N} \rightarrow \mathbb{R}^{N},
\end{equation*}
 and $f=f(x,t;  \mathbf{X} )$ is the solution to the Cauchy Problem in $\Omega=\mathbb{R}^{N} \times [0,T]$ 
\begin{equation}\label{peq}
\left\{
\begin{array}{rl}
Lf(x,t;  \mathbf{X} )&=g(x,\mathbf{X}(t)), \ \ \ \ \ \ \ \ (x,t) \in \mathbb{R}^N \times (0,T),\\
f(x,0)&= \varphi(x),  \ \ \ \ x \in \mathbb{R}^N,
\end{array}
\right.
\end{equation}
with
\begin{equation*}
\varphi: \mathbb{R}^N \longrightarrow  \mathbb{R} , \ \ g: \mathbb{R}^N \times \mathbb{R}^{N \times n} \longrightarrow  \mathbb{R} 
 \end{equation*}  continuous functions, 
and $L$ the following differential parabolic-type operator:

\begin{equation}\label{Lparabolic}
L = \displaystyle \sum_{i,j=1}^{N} a_{i,j}(x,t) \partial_{x_i,x_j}^{2}+ 
	 \sum_{i=1}^{N} b_{i}(x,t) \partial_{x_i}+c(x,t)- \partial_{t}.
\end{equation}
Through out all the paper, we will denote \mbox{$\mathbf{X}(t):=\left[\mathbf{x}_1(t),...,\mathbf{x}_n(t)\right]$}, $ \mathbf{V}(t):=\left[\mathbf{v}_1(t),...,\mathbf{v}_n(t)\right]$ $\in \mathbb{R}^{N \times n}$, 
and \mbox{$\mathbf{F}:=\left[\mathbf{F}_1,...,\mathbf{F}_n\right]$}.
In particular, we will denote $\mathbf{X}(0)=\mathbf{X}_0$ and $\mathbf{V}(0)=\mathbf{V}_0$.\\
We shall consider the following standing assumptions:

\begin{itemize}

\item[H1)] For every $i=1,...,n$, $\mathbf{F}_i$ 
is a continuous function, satisfying the following condition:
$\exists L_F \ge 0 $ such that $\forall K \subset \mathbb{R}^{N\times n}$ compact, $\exists L_{F}^{K} \ge 0 $ such that 
\begin{equation}\label{Flip}
\left.
\begin{array}{rl}
\left| \mathbf{F}_{i}(t,\mathbf{X},\mathbf{V}, \mathbf{w})-\mathbf{F}_{i}(t,\mathbf{\widehat{X}},\mathbf{\widehat{V}}, \mathbf{\widehat{w}})\right| \le L_{F}^{K} \left( \left| \mathbf{X}-\mathbf{\widehat{X}} \right| + \left| \mathbf{V}-\mathbf{\widehat{V}}\right|\right)+L_{F} \left| \mathbf{w}-\mathbf{\widehat{w}}  \right|,
 \end{array}
 \right.
\end{equation}
for any $\mathbf{X}, \mathbf{\widehat{X}}, \mathbf{V}, \mathbf{\widehat{V}} \in K$, $\mathbf{w}, \mathbf{\widehat{w}} \in \mathbb{R}^{N} $, $t \in \left[0,T\right]$.
\item[H2)]
 
The coefficients $a_{i,j}$, $b_{i}$, $c$, are bounded H\"older continuous function in $\Omega$, with coefficient $\alpha \in (0,1)$ with respect to $x$ and $\alpha/2$ with respect to $t$.\\
\item[H3)]  $L$ is uniformly parabolic in $\Omega$, meaning that
there exist $\mu_0,\mu_1>0$ such that, for every $\xi \in \mathbb{R}^{N}$ it holds
\begin{equation}\label{unifpar}
\mu_0 \left | \xi \right |^2 \le \sum_{i,j=1}^{N} a_{ij}(x,t) \xi_i \xi_j
\le \mu_1 \left | \xi \right |^2 \ \ \ \ \forall (x,t) \in \Omega.
\end{equation}
Moreover, there exists a constant $0\leq C<\lambda_0/4T$, with $\lambda_0\leq \mu_0/\mu_1^2$ (See  \cite{friedman} and also Appendix for more details) such that:\\
\item[H4)]
\begin{equation}\label{condH4}
|\varphi(x)-\varphi(\widehat{x})|\leq H \exp\left[C\max\left(|x|^2,|\widehat{x}|^2\right)\right] |x-\widehat{x}|^\alpha,
\end{equation}
for any $x$, $\widehat{x}\in\mathbb{R}^N$, for some constant $H\geq 0$.\\

\noindent

\item[H5)]

For every $R>0$, there exists a constant $H_R\geq 0$ such that
$$
\Big{|}g(x,{\bf X})-g(\widehat{x},\widehat{\bf X})\Big{|}\leq H_R \exp\left[C\max\left(|x|^2,|\widehat{x}|^2\right)\right]\left\{|x-\widehat{x}|^\alpha+\Big{|}{\bf X}-\widehat{\bf X}\Big{|}\right\}
$$
for any $x$, $\widehat{x}\in\mathbb{R}^N$, ${\bf X}$, $\widehat{\bf X}\in B_R$.
\end{itemize}
In the following we shall make use of the following amount $\ell(\theta,\nu):=e^{-\theta} (\theta/\nu)^\theta$ for every $\theta>0$, $\nu>0$ that is easy to check that corresponds to the global maximum of the function $y\in [0,\infty) \mapsto y^{\theta} e^{-\nu y}$. \\

 \noindent
\textbf{Remark 1.}
If Assumption H5) holds true, for every ${\bf X}={\bf X}(t)$ $\mathbb{R}^{N\times n}$-valued continuous function on $[0,T]$, we define $g_{\bf X}(x,t)$ $:=g(x,{\bf X}(t))$, which has the following properties: 
let $\bar{R}=\sup_{t\in [0,T]} |{\bf X}(t)|$, then, for every $M>0$,
\begin{equation}\label{growthg}
|g_{\bf X}(x,t)-g_{\bf X}(\widehat{x},t)|\leq H_{\bar{R},M} |x-\widehat{x}|^{\alpha},
\end{equation}
for any $x$, $\widehat{x}\in\mathbb{R}^{N}$ such that $|x|$, $|\widehat{x}|\leq M$, $t\in [0,T]$,
where $H_{\bar{R},M}=H_{\bar{R}} \exp[CM^2]$, that is $g_{\bf X}$ is locally H\"older continuous in $x$ with exponent $\alpha$, uniformly with respect to $t$.\\
\noindent Moreover, let $C< C^{\prime}<\lambda_0/4T$, then, from H5), we obtain
\begin{align}
& |g_{\bf X}(x,t)|\leq |g(0,0)|+H_{\bar {R}} e^{C|x|^2}\left(|x|^{\alpha}+{\bar R}\right)\nonumber\\
&\leq |g(0,0)|e^{C^\prime|x|^2}+H_{\bar{R}} e^{C^{\prime}|x|^2}\left(|x|^{\alpha}e^{-(C^\prime-C)|x|^2}+\bar{R} e^{C^{\prime}|x|^2}\right)\nonumber\\
&\leq \left[|g(0,0)|+H_{\bar {R}} (\ell(\alpha/2,C^\prime-C)+{\bar R})\right] e^{C^{\prime}|x|^2}. 
\end{align}
Similarly, such an upper bound is satisfied by the function $\varphi$, thanks to the Assumption H4) . These properties, satisfied by $\varphi$ and $g_{\bf X}$, in light of Theorem 12 (page 25) in \cite{friedman}, allow to establish that $f\left(\cdot, \cdot; \mathbf{X}\right): \mathbb{R}^{N} \times \left[ 0,T\right] \rightarrow \mathbb{R}$ defined by
\begin{equation}\label{fgamma}
f\left(x,t; \mathbf{X}\right)= \int_{\mathbb{R}^{N}} \Gamma \left( x,t;\xi,0\right)\varphi(\xi) d\xi - \int_{0}^{t}\int_{\mathbb{R}^{N}} \Gamma \left( x,t;\xi,\tau\right)g_{\bf X}(\xi,\tau) d\xi d\tau
\end{equation}
is the unique solution to the Cauchy problem ($\ref{peq}$) associated to $ \mathbf{X}$, meaning in particular that $f$ is continuous on $\mathbb{R}^N \times \left[0,T \right]$, $\partial_t f$,$\partial_{x_i} f$,$\partial_{x_i,x_j} f$ are continuous on $\mathbb{R}^N \times \left(0,T \right)$. Here the function $\Gamma$ is a fundamental solution of $Lu=0$ (See \cite{friedman} and Appendix for the details). \\
Later on, we will use estimates concerning $\Gamma$ and its derivatives:\\
$\forall \lambda_0^{*}<\lambda_0$ there exists a constant $C_{\Gamma}$ such that
\begin{equation}\label{solo_gamma}
\left| \displaystyle \Gamma  \left( x,t;\xi,\tau\right) \right| \le C_\Gamma \frac{1}{\left(t-\tau\right)^{\frac{N}{2}}} \displaystyle e^{-\frac{\lambda_0^{*}\left|x-\xi\right|^{2}}{4\left(t-\tau\right)}} \ \ \ \ \ \ \ \ \ \ \ i=1,...,N,
\end{equation}
\begin{equation}\label{der_gamma}
\left| \displaystyle \frac{\partial \Gamma \left( x,t;\xi,\tau\right)}{\partial x_i}\right| \le C_\Gamma \frac{1}{\left(t-\tau\right)^{\frac{N+1}{2}}} \displaystyle e^{-\frac{\lambda_0^{*}\left|x-\xi\right|^{2}}{4\left(t-\tau\right)}} \ \ \ \ \ \ \ i=1,...,N,
\end{equation}
\begin{equation}\label{der_seconda_gamma}
\left| \displaystyle \frac{\partial \Gamma \left( x,t;\xi,\tau\right)}{\partial x_i x_j}\right| \le C_\Gamma \frac{1}{\left(t-\tau\right)^{\frac{N+2}{2}}} \displaystyle e^{-\frac{\lambda_0^{*}\left|x-\xi\right|^{2}}{4\left(t-\tau\right)}}
 \ \ \ \ \ \ i,j=1,...,N.
\end{equation}
for any $x, \xi \in \mathbb{R}^{N}$, $0 \le \tau < t \le T$.\\
For (\ref{solo_gamma})-(\ref{der_gamma}) the reader is referred to \cite{friedman2}, Theorem 4.5, whereas
for (\ref{der_seconda_gamma}) to  \cite{pascucci2}.
We observe that in \cite{pascucci2}, inequality (\ref{der_seconda_gamma}) is proved for a hypoelliptic differential operator.  
However, such estimate is clearly satisfied by the uniformly parabolic operator $L$ defined in (\ref{Lparabolic}).
An essential ingredient for subsequent results is the estimation of the first and second-order derivatives of $f(x,t;{\bf X})$ with respect to $x$, in terms of the supremum of ${\bf X}$ over $[0,T]$. The inequalities in Proposition $\ref{disug}$ improve the usual estimates available for the
solution of the Cauchy problem ($\ref{peq}$), thanks to the local H\"older continuity of the data. Here we get a precise estimate of the constants for such bounds.
\begin{Proposition} \label{disug}
Let the assumptions H2)-H5) be satisfied. Then, for every 
$0<\nu_0<\lambda_0/4-CT$, $x\in \mathbb{R}^N$, $0<t\leq T$, ${\bf X}\in C([0,T];\mathbb{R}^{N\times n})$, $i$, $j=1,\ldots, N$, the following inequalities hold true:
\begin{align}\label{Prop.main.1}
|\partial_{x_i} f(x,t;{\bf X})|\leq  Ke^{\kappa |x|^2}\left(\frac{H}{t^{\frac{1-\alpha}{2}}}+\frac{2}{\alpha+1} t^{\frac{\alpha+1}{2}} H_{\bf X}\right),
\end{align}
\begin{align}\label{Prop.main.2}
|\partial^2_{x_i x_j} f(x,t;{\bf X})|\leq  Ke^{\kappa |x|^2} \left(\frac{H}{t^{1-\frac{\alpha}{2}}}+\frac{2}{\alpha} t^{\frac{\alpha}{2}} H_{\bf X}\right),
\end{align}
where $C_{\Gamma}$ is obtained from inequalities (\ref{solo_gamma})-(\ref{der_seconda_gamma}) for $\lambda_0^*=\lambda_0-4\nu_0$, $H_{\mathbf{X}}$ stands for the constant in $H5$), for $R=\|X\|_{\infty,T}$, $H$ the one in H4), and
\begin{align}\label{Prop.main.constants}
K:=\frac{\pi^{N/2}C_{\Gamma}\ell(\alpha/2,\nu_0)}{\left[\lambda_0/2-2CT\right]^{N/2}}, \qquad\qquad \kappa:=\frac{C^2 T}{\lambda_0/4-CT}+2C.
\end{align}
\end{Proposition}

\noindent
{\em Proof.} We focus on (\ref{Prop.main.1}). First we recall that the Fundamental solution $\Gamma$ satisfies
\begin{align}\label{Prop.1}
\int_{\mathbb{R}^N} \Gamma(x,t,\xi,\tau)d\xi=1,
\end{align}
for any $x\in \mathbb{R}^N$, $0\leq \tau<t\leq T$. Therefore, we have
\begin{equation}\label{Prop.2}
0=\partial_{x_i} \int_{\mathbb{R}^N}\Gamma(x,t,\xi,\tau)d\xi=\int_{\mathbb{R}^N}\partial_{x_i}\Gamma(x,t,\xi,\tau)d\xi,
\end{equation}
\begin{equation}\label{Prop.22}
0=\partial^2_{x_i x_j} \int_{\mathbb{R}^N}\Gamma(x,t,\xi,\tau)d\xi=\int_{\mathbb{R}^N}\partial^2_{x_i x_{j}}\Gamma(x,t,\xi,\tau)d\xi.
\end{equation}
Let ${\bf X}\in C([0,T];\mathbb{R}^{N\times n})$ and set $\overline{R}=\|{\bf X}\|_{\infty,T}$. Then, the representation formula (\ref{fgamma}) implies, for every $x\in \mathbb{R}^N$, $t\in (0,T]$, the inequality $|\partial_{x_i} f(x,t;{\bf X})|\leq G_1+G_2$, where 
\begin{align}\label{Prop.3}
G_1:=\Big{|}\int_{\mathbb{R}^N} \partial_{x_i} \Gamma(x,t,\xi,0)\varphi(\xi) d\xi-\varphi(x)\int_{\mathbb{R}^N} \partial_{x_i} \Gamma(x,t,\xi,0) d\xi \Big{|},
\end{align}
\begin{align}\label{Prop.4}
G_2:=\Big{|}\int_0^t\int_{\mathbb{R}^N} \partial_{x_i} \Gamma(x,t,\xi,\tau)g_{\bf X}(\xi,\tau) d\xi-\int_0^t g_{\bf X}(x,\tau)\int_{\mathbb{R}^N} \partial_{x_i} \Gamma(x,t,\xi,\tau) d\xi d\tau\Big{|}.
\end{align}
By the inequalities in $H4$)-$H5$), since $C<\lambda_0/4T$, we choose $\lambda_0^*\in (4TC, \lambda_0)$ and we set $\nu_0=(\lambda_0-\lambda_0^*)/4$. Therefore from (\ref{der_gamma}), for $0<t\leq T$, we get
\begin{align}\label{Prop.5}
G_1 &\leq \frac{C_{\Gamma}H}{t^\frac{N+1}{2}}\int_{\mathbb{R}^N} e^{-\frac{\lambda_0^*}{4}\frac{|x-\xi |^2}{t}+C|x|^2+C|\xi|^2} |x-\xi|^{\alpha} d\xi \nonumber\\
&\leq C_{\Gamma}H\, t^{\frac{\alpha-1}{2}} \int_{\mathbb{R}^N} e^{-\left(\lambda_0^*/4-Ct\right)|u|^2+2C\sqrt{t}\langle x,u\rangle+2C|x|^2} |u|^{\alpha} du\nonumber\\
&\leq C_{\Gamma}H \ell(\alpha/2,\nu_0)\, t^{\frac{\alpha-1}{2}} \int_{\mathbb{R}^N} 
e^{-\left(\lambda_0^*/4-Ct+\nu_0\right)|u|^2+2C\sqrt{t}\langle x,u\rangle+2C|x|^2}  du\nonumber\\
& \leq C_{\Gamma} H \ell(\alpha/2,\nu_0) t^{\frac{\alpha-1}{2}} e^{\frac{C^2 t |x|^2 }{\lambda_0/4-Ct}}\int_{\mathbb{R}^N} e^{-\Big{|}u \sqrt{\lambda_0/4-Ct}-\frac{C\sqrt{t} x}{\sqrt{\lambda_0/4-Ct}}\Big{|}^2+2C|x|^2} du \nonumber\\
& =\frac{C_{\Gamma}H\ell(\alpha/2,\nu_0)}{t^{\frac{1-\alpha}{2}}\left[\lambda_0/4-Ct\right]^{N/2}} e^{\left(\frac{C^2 t}{\lambda_0/4-Ct}+2C\right)|x|^2} \left(\frac{\pi}{2}\right)^{N/2}.
\end{align}
Here we have applied the changes of variable $\xi =x+\sqrt{t} u$ and $v/2=u \sqrt{\lambda_0/4-Ct}$ $-\frac{C\sqrt{t} x}{\sqrt{\lambda_0/4-Ct}}$ and we have used the well known relation $\int_{\mathbb{R}^N} e^{-\frac{1}{2}|v^2|} dv=[2\pi]^{N/2}$. Since $t\in \mapsto [\lambda_0/4-Ct]^{-1}$ is an increasing function, we finally obtain the estimate
\begin{align}\label{Prop.6}
G_1 & \leq \frac{C_{\Gamma}H\ell(\alpha/2,\nu_0)}{t^{\frac{1-\alpha}{2}}\left[\lambda_0/4-CT\right]^{N/2}} e^{\left(\frac{C^2 T}{\lambda_0/4-CT}+2C\right)|x|^2} \left(\frac{\pi}{2}\right)^{N/2}.
\end{align}
By similar arguments applied to $G_2$, using $H5$) (with ${\bf X}\equiv\widehat{\bf X}$) and replacing $t$ with $t-\tau$ in the integral over $\mathbb{R}^N$, we can write
\begin{align}\label{Prop.7}
G_2 &\leq \frac{C_{\Gamma}H_{\overline{R}}\ell(\alpha/2,\nu_0)}{\left[\lambda_0/4-CT\right]^{N/2}} e^{\left(\frac{C^2 T}{\lambda_0/4-CT}+2C\right)|x|^2} \left(\frac{\pi}{2}\right)^{N/2} \int_0^t \frac{1}{(t-\tau)^{\frac{1-\alpha}{2}}} d\tau\nonumber\\
&= \frac{2C_{\Gamma}H_{\overline{R}}\ell(\alpha/2,\nu_0)}{(\alpha+1)\left[\lambda_0/4-CT\right]^{N/2}} e^{\left(\frac{C^2 T}{\lambda_0/4-CT}+2C\right)|x|^2} \left(\frac{\pi}{2}\right)^{N/2}  t^{\frac{\alpha+1}{2}},
\end{align}
Inequalities (\ref{Prop.6}) and (\ref{Prop.7}) yield (\ref{Prop.main.1}). We observe that the proof of inequality (\ref{Prop.main.2}) follows through similar passages using the estimate for the second-order derivatives of the fundamental solution in (\ref{der_seconda_gamma}) together with (\ref{Prop.22}), hence we omit the details.\\


\section{Existence and Uniqueness of a local solution.}

Using a fixed point argument, we prove the local existence and uniqueness of the solution of $(\ref{dinamica})$-$(\ref{peq})$. We state more clearly that, by solution of $(\ref{dinamica})$-$(\ref{peq})$, we refer to the couple $\mathbf{Y}=\left(\mathbf{X},\mathbf{V}\right)$, where $\mathbf{X},\mathbf{V} \in C\left([0,T];\mathbb{R}^{N\times n}\right) \cap C^{1}\left((0,T);\mathbb{R}^{N\times n}\right) $, and $f \in C^{2,1}\left(\mathbb{R}^{N} \times (0,T)\right)$, satisfying $(\ref{peq})$, is expressed as in $(\ref{fgamma})$.
\begin{Theorem}\label{teolocale}
Under hypotheses H1)-H5), system $(\ref{dinamica})$-$(\ref{peq})$ has a unique solution on $[0,\overline{T}]$, where $\overline{T} \in \left(0,T\right]$ depends on $\mathbf{X_0}$, $\mathbf{V_0}$, $\alpha$, $n$, $N$, $R$.
\end{Theorem}
%

%
%
%
\noindent
{\em Proof.} Let $R>0$ and $0<\overline{T}\le T$.\\
In the following we denote \mbox{$C_0:=\displaystyle \max_{\tau \in \left[0,\overline{T}\right] }  | \mathbf{F}(\tau,\mathbf{X}_0,\mathbf{V}_0, 0) |$}. 
We consider the mapping $\Psi$ defined as follows

\begin{equation}\label{operatore}
\Psi\left(\mathbf{X}, \mathbf{V}\right) (t) = 
\begin{pmatrix}
{\mathbf{x}_0}_1+ \displaystyle \int_{0}^{t}\mathbf{v}_1 (\tau) \, d \tau \\
\vdots \\
{\mathbf{x}_0}_n+ \displaystyle \int_{0}^{t}\mathbf{v}_n (\tau) \, d \tau \\
{\mathbf{v}_0}_1+ \displaystyle \int_{0}^{t}\mathbf{F}_1 \left( \tau, \mathbf{X}(\tau), \mathbf{V}(\tau), \nabla f\left( \mathbf{x}_1(\tau), \tau; \mathbf{X} \right)\right) \, d \tau \\
\vdots \\
{\mathbf{v}_0}_n+ \displaystyle \int_{0}^{t}\mathbf{F}_n \left( \tau, \mathbf{X}(\tau), \mathbf{V}(\tau), \nabla f\left( \mathbf{x}_n(\tau), \tau; \mathbf{X} \right)\right) \, d \tau \\
\end{pmatrix}
\end{equation}
for any $ \left( \mathbf{X}, \mathbf{V}\right) \in E_R $, $t \in [0,\overline{T}]$, where $E_R:=C\left( [0,\overline{T}];  B_{R}(\mathbf{X}_0) \times B_{R}(\mathbf{V}_0) \right)$, and $f$ is the function in (\ref{fgamma}).\\
Due to (\ref{Flip}) and (\ref{Prop.main.1}), we have that $\Psi\left(\mathbf{X}, \mathbf{V}\right)$ is continuous at $t=0$, thus 
$$
\Psi\left(\mathbf{X}, \mathbf{V}\right) \in C\left( [0,\overline{T}]; \mathbb{R}^{N\times n} \times \mathbb{R}^{N\times n} \right).
$$
We observe that a fixed point of $\Psi$, $\left( \mathbf{\overline{X}}, \mathbf{\overline{V}} \right) \in E_R$, is a solution to (\ref{dinamica})-(\ref{peq}).\\
We shall define suitable conditions on $\overline{T}$, in order to guarantee 
$
\Psi: E_R \longrightarrow E_R,
$
and the fact that $\Psi$ is a contraction operator. Since $E_R$, endowed with the uniform norm, is clearly a Banach space, the result follows from fixed-point theorem (\cite{evans}, Theorem 1, p.534).


Let $ \left(\mathbf{X},  \mathbf{V}\right) \in E_R$. We consider the first $n$ components of the operator $\Psi$. Since we assume $\mathbf{V} \in B_{R}(\mathbf{V}_0) $, it follows that 
$$
\left | \Psi_j\left(\mathbf{X}, \mathbf{V}\right) (t)- {\mathbf{x}_0}_j \right | \le \overline{T} \left( R + | \mathbf{V}_0| \right) \ \ \ \ \forall t \in [0, \overline{T}], j=1,...,n.
$$
Thus
\begin{equation}\label{primecomp}
\left | 
\left(
\begin{array}{lll}
\Psi_1\left(\mathbf{X}, \mathbf{V}\right)(t)- {\mathbf{x}_0}_{1}\\
\vdots \\
\Psi_n\left(\mathbf{X}, \mathbf{V}\right)(t)-  {\mathbf{x}_0}_{n}
\end{array}
\right)
 \right |
 \le n \overline{T} \left( R + | \mathbf{V}_0| \right).
\end{equation}
Hence $\left[\Psi_1\left(\mathbf{X}, \mathbf{V}\right)(t),...,\Psi_n\left(\mathbf{X}, \mathbf{V}\right)(t)\right] \in B_R\left(\mathbf{X}_0\right)$ for any $t \in [0, \overline{T}]$ if we impose the condition:
\begin{equation}\label{cond1T}
\overline{T} \le \displaystyle \frac{R}{n (R + | \mathbf{V}_0|)}.
\end{equation}
Let us consider the component $l \in \left\{n+1,...,2n\right\}$.
Recalling ($\ref{Flip}$) and ($\ref{Prop.main.1}$) , we obtain in particular that 
\begin{equation}\label{comp_prov}
\left.
\begin{array}{rl}
&\left | \Psi_{j+n}\left(\mathbf{X}, \mathbf{V}\right) (t)- {\mathbf{v}_0}_{j} \right | \le
\displaystyle \int_{0}^{t} \left| \mathbf{F}_{j}(\tau,\mathbf{X}(\tau),\mathbf{V}(\tau), \nabla f (\mathbf{x}_{j}(\tau), \tau;\mathbf{X} ))\right| d\tau \le\\\\
&\le \displaystyle \int_{0}^{t} \left| \mathbf{F}_{j}(\tau,\mathbf{X}(\tau),\mathbf{V}(\tau), \nabla f (\mathbf{x}_{j}(\tau), \tau;\mathbf{X} ))- \mathbf{F}_{j}(\tau,\mathbf{X}_0,\mathbf{V}_0, 0)\right| + | \mathbf{F}_{j}(\tau,\mathbf{X}_0,\mathbf{V}_0, 0) | d\tau \le\\\\
&\le 2 L_{F}^{R} R\overline{T}+ \displaystyle \int_{0}^{t}L_{F} \left |\nabla f (\mathbf{x}_{j}(\tau), \tau;\mathbf{X} ) \right|+|\mathbf{F}_{j}(\tau,\mathbf{X}_0,\mathbf{V}_0, 0) |  d\tau \le\\\\
&\le \displaystyle 2 L_{F}^{R} R\overline{T}+ \displaystyle \int_{0}^{t} L_{F} \sqrt{N}Ke^{\kappa \left( |\mathbf{X}_0|^{2}+R^{2}\right)}\left(\frac{H}{\tau^{\frac{1-\alpha}{2}}}+\frac{2}{\alpha+1} \tau^{\frac{\alpha+1}{2}} H_{\bf X}\right) + | \mathbf{F}_{j}(\tau, \mathbf{X}_0,\mathbf{V}_0, 0) | d \tau \le \\\\
& \le \overline{T} \left(2 L_{F}^{R} R +C_0  \right) + L_{F}\sqrt{N}Ke^{\kappa \left( |\mathbf{X}_0|^{2}+R^{2}\right)} 
\left( \displaystyle \frac{2H}{\alpha+1}\overline{T}^{\frac{\alpha+1}{2}}+ \frac{2 H_{\bf X}}{\alpha+1} \frac{2}{\alpha+3}\overline{T}^{\frac{\alpha+3}{2}}\right).
\end{array}
\right.
\end{equation}
Thus
\begin{equation*}
\left | 
\left(
\begin{array}{lll}
\Psi_{n+1}\left(\mathbf{X}, \mathbf{V}\right)- {\mathbf{v}_0}_{1}\\
\vdots \\
\Psi_{2n}\left(\mathbf{X}, \mathbf{V}\right)-  {\mathbf{v}_0}_{n}
\end{array}
\right)
 \right | \le
 \end{equation*}

\begin{equation}\label{secondecomp}
\le n \overline{T} \left( 2 L_{F}^{R} R + C_0  \right) + nL_{F}\sqrt{N}Ke^{\kappa \left( |\mathbf{X}_0|^{2}+R^{2}\right)} 
\left( \displaystyle \frac{2H}{\alpha+1}\overline{T}^{\frac{\alpha+1}{2}}+ \frac{2 H_{\bf X}}{\alpha+1} \frac{2}{\alpha+3}\overline{T} \overline{T}^{\frac{\alpha+1}{2}}\right).
\end{equation}
From $(\ref{cond1T})$ follows, in particular, that $\overline{T}<1$, thus $\overline{T}<\overline{T}^{\frac{\alpha+1}{2}} $ since $\alpha \in \left(0,1\right)$. Hence it sufficies to impose 
$\overline{T} \le T_1$,
%
where $T_1=T_1\left( \mathbf{X}_0, \mathbf{V}_0, \alpha, n, N, R\right)$ is defined as
\begin{equation}\label{defT1}
T_1 := \min \left( \displaystyle \frac{R}{n (R + | \mathbf{V}_0|)}, \displaystyle \left(\frac{R}{\frac{2n\sqrt{N}L_{F}Ke^{\kappa \left( |\mathbf{X}_0|^{2}+R^{2}\right)}}{\alpha+1}\left(1+\frac{2 H_{\bf X}}{\alpha +3}\right)} \right)^{\frac{2}{\alpha+1}}
  \right),
\end{equation}
to ensure that the range of $\Psi$ is a subset of E.\\
We now show that it is possible to determine a condition on $\overline{T}$ such that $\Psi$ is a contraction operator. 
Let $ \left(\mathbf{X}, \mathbf{V} \right)$ and $\left( \widehat{\mathbf{X}}, \widehat{\mathbf{V}} \right)  \in E$.\\
%
Clearly it holds 
\begin{equation}\label{contra_prime}
\left.
\begin{array}{rl}
&\left | \Psi_j\left(\mathbf{X}, \mathbf{V}\right) (t)- \Psi_j\left(\widehat{\mathbf{X}}, \widehat{\mathbf{V}}\right) (t) \right | \le
\displaystyle \int_{0}^{t} \left | \mathbf{v}_{j}(\tau) - \hat{\mathbf{v}}_{j}(\tau)  \right | d \tau\\\\
&\le \overline{T} ||  \mathbf{V} - \widehat{\mathbf{V}} ||_{\infty,\overline{T}} \le
 \overline{T} \left|\left| \left( \mathbf{X}, \mathbf{V}\right) - \left( \widehat{\mathbf{X}}, \widehat{\mathbf{V}}\right) \right|\right|_{\infty,\overline{T}}.
\end{array}
\right.
\end{equation}
We now focus on the remaining components. From the Lipschitz condition ($\ref{Flip}$), it follows that, for $j=1,...,n$: 
\begin{equation}\label{contra_second}
\left.
\begin{array}{rl}
&\left |  \Psi_{n+j}\left(\mathbf{X}, \mathbf{V}\right) (t)-\Psi_{n+j}\left(\widehat{\mathbf{X}}, \widehat{\mathbf{V}}\right) (t) \right | \le\\\\
&\le L_{F}^{R} \displaystyle \int_{0}^{t} \left |  \mathbf{X}(\tau)- \widehat{\mathbf{X}}(\tau)\right | + \left |  \mathbf{V}(\tau)- \widehat{\mathbf{V}}(\tau)\right | d\tau \\\\
&+ L_{F} \displaystyle \int_{0}^{t} \left | \nabla f (\mathbf{x}_{j}(\tau), \tau;\mathbf{X} ))- \nabla f (\widehat{\mathbf{x}}_{j}(\tau), \tau;\widehat{\mathbf{X}} ))\right| d\tau \le \\\\
& \le L_{F}^{R}  \overline{T} \left(  || \mathbf{X} - \widehat{\mathbf{X}}||_{\infty,\overline{T}}+ || \mathbf{V} - \widehat{\mathbf{V}}||_{\infty, \overline{T}}\right) + L_{F} \left( I_1 + I_2 \right),
\\\\
\end{array}
\right.
\end{equation}
where $ I_1, I_2$ are the integrals
\begin{equation}
I_1 :=  \displaystyle \int_{0}^{t} \left | \nabla f (\mathbf{x}_{j}(\tau), \tau;\mathbf{X} ))- \nabla f (\widehat{\mathbf{x}}_{j}(\tau), \tau;\mathbf{X} ))\right| d\tau,
\end{equation}
\begin{equation}
I_2 := \displaystyle \int_{0}^{t} \left | \nabla f (\widehat{\mathbf{x}}_{j}(\tau), \tau;\mathbf{X} ))- \nabla f (\widehat{\mathbf{x}}_{j}(\tau), \tau; \widehat{\mathbf{X}}))\right| d\tau.
\end{equation}
We observe that, by the mean value theorem and $(\ref{Prop.main.2})$ , it holds
\begin{equation}\label{valor_medio}
\left.
\begin{array}{rl}
&\left | \nabla f (\mathbf{x}_{j}(\tau), \tau;\mathbf{X} ))- \nabla f (\widehat{\mathbf{x}}_{j}(\tau), \tau;\mathbf{X} ))\right|  \le\\\\
& \le \displaystyle \sum_{i=1}^{N} \left|  \partial_{x_i}f (\mathbf{x}_{j}(\tau), \tau;\mathbf{X} ))-\partial_{x_i}f (\widehat{\mathbf{x}}_{j}(\tau), \tau;\mathbf{X} )) \right| \le\\\\
&  \le \displaystyle \sum_{i=1}^{N} \sum_{m=1}^{N} \left|  \partial_{x_{im}}^{2} f (\mathbf{x}_{ij}^{*}, \tau;\mathbf{X} )\right| \left( \left|\mathbf{x}_{j}(\tau)- \widehat{\mathbf{x}}_{j}(\tau)\right| \right) \le\\\\
&  \le \displaystyle N^2 Ke^{\kappa  \left | \mathbf{x}_{ij}^{*} \right |^{2}}\left( \frac{H}{\tau^{1-\alpha/2}}+ \frac{2 H_{\bf X}}{\alpha}\tau^{\alpha/2}\right) \left( \left|\mathbf{x}_{j}(\tau)- \widehat{\mathbf{x}}_{j}(\tau)\right| \right) ,
\end{array}
\right.
\end{equation}
where ${\mathbf{x}^{*}}_{ij}$  belongs to the segment connecting $\mathbf{x}_{j}(\tau)$ to $\widehat{\mathbf{x}}_{j}(\tau)$ thus, in particular, ${\mathbf{x}^{*}}_{ij} \in B_R\left( \mathbf{X}_0\right)$. From the parallelogram inequality we get $\left | {\mathbf{x}^{*}}_{ij} \right |^{2} \le 2 \left(|\mathbf{X}_0|^{2}+R^{2}\right).$
We finally achieve the following estimate for $I_1$:
\begin{equation}\label{stimI1}
I_1 \le N^2 Ke^{\kappa \left( \left | \mathbf{X}_0 \right |^{2}+R^2\right)}\left( \frac{2H}{\alpha} \overline{T}^{\alpha/2}+ \frac{2 H_{\bf X}}{\alpha}\overline{T}^{\alpha/2+1}\right)  \left|\left| \mathbf{X} - \widehat{\mathbf{X}}\right|\right|_{\infty,\overline{T}}. 
\end{equation}
We now focus on $I_2$.
Since $f$ is a solution of ($\ref{peq}$), we observe that $\eta(\xi,\tau):=f(\xi,\tau;\mathbf{X})-f(\xi,\tau;\widehat{\mathbf{X}})$ satisfies the following Cauchy problem:
\begin{equation}\label{peqdiff}
\left\{
\begin{array}{rl}
L \eta (\xi, \tau)&=g (\xi, \mathbf{X}(\tau))-g (\xi, \widehat{\mathbf{X}}(\tau)) \ \ \ \ \ \ \ \ (\xi,\tau) \in \mathbb{R}^N \times (0,T),\\
\eta(\xi,0)&= 0 \ \ \ \ \ \ \ \ \ \ \ \ \ \ \ \ \ \ \ \ \ \ \ \ \ \ \ \ \ \ \ \ \ \ \ \ \ \ \xi \in \mathbb{R}^N.
\end{array}
\right.
\end{equation}
$I_2$ can thus be rewritten as $I_2= \displaystyle \int_{0}^{t} \left| \nabla \eta (\widehat{\mathbf{x}}_{j}(\tau),\tau)\right| d \tau$.\\
From ($\ref{der_gamma}$) and H5) we have
\begin{equation}\label{stim1}
\left.
\begin{array}{rl}
&\left| \nabla \eta (\xi,\tau)\right| =
 \displaystyle \left| \int_{0}^{\tau} \int_{\mathbb{R}^{N}} \nabla \Gamma (\xi,\tau; \overline{\xi},\overline{\tau}) \{g(\overline{\xi},\mathbf{X}(\tau))-g(\overline{\xi},\widehat{\mathbf{X}}(\tau))\} d \overline{\xi} d\overline{\tau} \right|\le \\\\
 & \le C_{\Gamma}\displaystyle \int_{0}^{\tau} \int_{\mathbb{R}^{N}} \frac{e^{-\frac{\lambda_{0}^{*}}{4}\frac{\left| \xi -\overline{\xi} \right|^{2}}{\tau-\overline{\tau}}}}{\left(\tau - \overline{\tau}\right)^{\frac{N+1}{2}}}  \left | g(\overline{\xi},\mathbf{X}(\tau))-g(\overline{\xi},\widehat{\mathbf{X}}(\tau)) \right|d \overline{\xi} d\overline{\tau} \le \\\\
 & \le C_{\Gamma}H_R \displaystyle \int_{0}^{\tau} \int_{\mathbb{R}^{N}}  \frac{e^{-\frac{\lambda_{0}^{*}}{4}\frac{\left| \xi -\overline{\xi} \right|^{2}}{\tau-\overline{\tau}}}}{\left(\tau - \overline{\tau}\right)^{\frac{N+1}{2}}}  e^{C \left | \overline{\xi} \right |^{2}} d \overline{\xi} d\overline{\tau} \left | \left| \mathbf{X}- \widehat{\mathbf{X}} \right | \right |_{\infty,\overline{T}} = \\\\
 & = C_{\Gamma} H_{R} \displaystyle \int_{0}^{\tau} \int_{\mathbb{R}^{N}}  \frac{e^{-\frac{\lambda_{0}^{*}}{4}\left| y \right|^{2}+C\left | \xi + \sqrt{\tau-\overline{\tau}} y \right |^{2}}}{\left(\tau - \overline{\tau}\right)^{\frac{N+1}{2}}}  \left( \tau-\overline{\tau}\right)^{\frac{N}{2}}dy d\overline{\tau}\left | \left| \mathbf{X}- \widehat{\mathbf{X}} \right | \right |_{\infty,\overline{T}}
 \end{array}
\right.
\end{equation}
where, in the last equality, we have performed the change of variable  
\begin{equation}\label{cambio_coord1}
 \displaystyle y=\frac{-\xi+\overline{\xi}}{ \sqrt{\tau-\overline{\tau}}} .
 \end{equation}
Since  $\left | \xi + \sqrt{\tau-\overline{\tau}} y \right |^{2} \le 2 \left| \xi \right |^{2} + 2 \left| y \right |^{2}$,
we rewrite the exponential term in ($\ref{stim1}$). 
We require $ \displaystyle \overline{T}< \frac{\lambda_{0}^{*}}{8C}$, implying that $\overline{\gamma}:=\left (\frac{\lambda_0^{*}}{4}-2C\overline{T}\right) >0$. Recalling the expression of $I_{0}\left(\overline{\gamma}\right)$, and its exact value (see Appendix for details), we get:
\begin{equation}\label{stim2}
\left.
\begin{array}{rl}
\left| \nabla \eta (\xi,\tau)\right| & \le C_{\Gamma} H_{R} e^{2 C\left| \xi \right|^{2}}  \displaystyle \int_{0}^{\tau} \int_{\mathbb{R}^{N}}  \frac{1}{\left(\tau - \overline{\tau}\right)^{\frac{1}{2}}} e^{- \overline{\gamma}\left| y \right|^{2}} dy d\overline{\tau} \left | \left |  \mathbf{X}- \widehat{\mathbf{X}} \right|\right |_{\infty,\overline{T}}=\\\\
&= C_{\Gamma} H_{R} e^{2 C\left| \xi \right|^{2}}\sqrt{\tau} \displaystyle \left(\frac{\pi}{\overline{\gamma}}\right)^{N/2} \left | \left |  \mathbf{X}- \widehat{\mathbf{X}} \right|\right |_{\infty,\overline{T}},
\end{array}
\right.
\end{equation}
By ($\ref{stim2}$), with $\left| \xi \right|^2=\left| \mathbf{x}_j(\tau) \right|^2 \le 2 \left| \mathbf{X}_0\right|^2+2R^2$, the estimate of $I_2$ immediately follows:
\begin{equation}\label{stim3}
I_{2} \le  C_{\Gamma} H_{R} e^{2C\left(\left|  \mathbf{X}_0\right|^2+R^2\right)}   \displaystyle \left(\frac{\pi}{\overline{\gamma}}\right)^{N/2} \overline{T}^{\frac{3}{2}}\left | \left |  \mathbf{X}- \widehat{\mathbf{X}} \right|\right |_{\infty,\overline{T}}.
\end{equation}
From $(\ref{contra_second})$, $(\ref{stimI1})$ and $(\ref{stim3})$ we get the following inequality 
\begin{equation}\label{contra_seconde}
\left.
\begin{array}{rl}
&\left |  \Psi_{n+j}\left(\mathbf{X}, \mathbf{V}\right)(t)  - \Psi_{n+j}\left(\widehat{\mathbf{X}}, \widehat{\mathbf{V}}\right) (t) \right | \le\\\\
& \le \left( 2L_{F}^{R} \overline{T} +L_{F} N^2 Ke^{\kappa \left( \left | \mathbf{X}_0 \right |^{2}+R^2\right)}\left( \frac{2H}{\alpha} \overline{T}^{\alpha/2}+ \frac{2 H_{\bf X}}{\alpha}\overline{T}^{\alpha/2+1}\right) +\right.\\\\
&+\left. L_{F} C_{\Gamma} H_{R} e^{2C\left(\left|  \mathbf{X}_0\right|^2+R^2\right)}  \displaystyle  \left(\frac{\pi}{\overline{\gamma}}\right)^{N/2} \overline{T}^{\frac{3}{2}} \right) \left | \left |  \left(\mathbf{X},\mathbf{V}\right)-\left( \widehat{\mathbf{X}},\widehat{\mathbf{V}}\right) \right|\right |_{\infty,\overline{T}},
\end{array}
\right.
\end{equation}
for any $\left(\mathbf{X},\mathbf{V}\right), \left( \widehat{\mathbf{X}},\widehat{\mathbf{V}}\right) \in E$, $j=1,...,n$, $t \in \left[0,T\right]$.\\
From ($\ref{contra_prime}$) and ($\ref{contra_seconde}$) we conclude that
\begin{equation}\label{contra_final}
\left.
\begin{array}{rl}
&\left | \left |  \Psi\left(\mathbf{X}, \mathbf{V}\right)  - \Psi \left(\widehat{\mathbf{X}}, \widehat{\mathbf{V}}\right) \right | \right |_{\infty,\overline{T}}\le\\\\
& \le \sqrt{2n} \left[ \max \left( 2L_{F}^{R},1 \right)\overline{T}+ S\left( \mathbf{X}_0, \mathbf{V}_0, \alpha, n, N, R, \overline{T}\right)  \right] \left | \left |  \left(\mathbf{X},\mathbf{V}\right)-\left( \widehat{\mathbf{X}},\widehat{\mathbf{V}}\right) \right|\right |_{\infty,\overline{T}},
\end{array}
\right.
\end{equation}
where 
\begin{equation}\label{defT2}
\left.
\begin{array}{rl}
& S \left( \mathbf{X}_0, \mathbf{V}_0, \alpha, n, N, R, \overline{T}\right) := \left( 2L_{F}^{R} \overline{T} +L_{F} N^2 Ke^{\kappa \left( \left | \mathbf{X}_0 \right |^{2}+R^2\right)}  \displaystyle \overline{T}^{\alpha/2}\frac{2}{\alpha}\left( H+H_{\bf X}\overline{T}\right)\right.\\
&+\left.L_{F} C_{\Gamma} H_{R} e^{2C\left(\left|  \mathbf{X}_0\right|^2+R^2\right)}  \displaystyle  \left(\frac{\pi}{\overline{\gamma}}\right)^{N/2} \overline{T}^{\frac{3}{2}} \right).
\end{array}
\right.
\end{equation}
Finally, a sufficient condition on $\overline{T}$ such that $\Psi$ is a contraction operator is given by
\begin{equation}\label{minT}
\overline{T} < \min \{T_1, T_2 \},
\end{equation}
 where $T_1= T_1\left( \mathbf{X}_0, \mathbf{V}_0, \alpha, n, N, R\right)$ is defined in ($\ref{defT1}$), and $T_2=T_2 \left( \mathbf{X}_0, \mathbf{V}_0, \alpha, n,\right.$ $\left. N, R\right)$ is such that $T_2  < \displaystyle \frac{\lambda_0^{*}}{8C}$ and $S \left( \mathbf{X}_0, \mathbf{V}_0, \alpha, n, N, R, T_2\right) <1$.
We observe that such a choice of $T_2$ can be obtained since $S$ is an increasing continuous function of $\overline{T}$, and $S \left( \mathbf{X}_0, \mathbf{V}_0, \alpha, n, N, R, 0\right)=0$.\\
Hence, thesis follows from fixed-point theorem (\cite{evans}, Theorem 1, p. 534), which ensures the existence and uniqueness of a fixed point for $\Psi$ over the interval $\left[0,\overline{T}\right]$.\\
%

\noindent
\textbf{Remark 2.} Let us consider the case of $\varphi \equiv 0$.
 From Proposition \ref{disug} we deduce that
\begin{equation}
\left| \nabla f \left( x,t; \mathbf{X}\right) \right|\le  \frac{2 \sqrt{N}}{\alpha+1} Ke^{\kappa |x|^2} T^{\frac{\alpha+1}{2}} H_{\bf X}.
\end{equation}
Therefore Assumption H1) can be weakened, requiring only the local Lipschitz continuity of $\mathbf{F}_i$ with respect to all the arguments.
Actually $\left | \nabla f \left(x,t;\mathbf{X} \right) \right| $ is bounded on $B_R \times \left[0, \overline{T}\right]\times E_R$, for any $R>0$.
\section{The case of a non-local concentration.}

 In the previous sections we established an existence result for a local solution to systems in the form ($\ref{dinamica}$)-($\ref{peq}$). Under the same assumptions, we consider a variant of system ($\ref{dinamica}$), given by
 \begin{equation}\label{dinamica_media}
\left\{
\begin{array}{rl}
\mathbf{\dot{x}}_i(t)&=\mathbf{v}_{i}(t), \ \ \ \ \mbox{$0<t<T$}\\\
\mathbf{\dot{v}}_i(t)&= \displaystyle \mathbf{F}_{i}\left( t, \mathbf{X}(t), \mathbf{V}(t), \fint_{B_{\delta}(\mathbf{x}_i(t))} \nabla f\left( \xi,t;\mathbf{X}\right) d\xi\right), \ \ \ \ \forall i=1,...,n
\end{array}
\right.
\end{equation}
where $f$ satisfies ($\ref{peq}$).
The main difference between ($\ref{dinamica}$) and ($\ref{dinamica_media}$) lies in the fact that  $ \nabla f \left( \mathbf{x}_i(t),t;  \mathbf{X} \right)$ is replaced by the average over a ball centered in $\mathbf{x}_i(t)$ and having radius $\delta >0$.
 A local existence result is stated in the following theorem.
 \begin{Theorem}
Under hypotheses H1)-H5), the system $(\ref{dinamica_media})$ has a unique solution on $[0,\widetilde{T}]$, where $\widetilde{T} \in \left(0,T\right)$ depends on $\mathbf{X_0}$, $\mathbf{V_0}$, $\alpha$, $n$, $N$, $R$.
 \end{Theorem}
 \noindent
{\em Proof.}  We consider the operator defined in $(\ref{operatore})$, replacing $\nabla f \left( \mathbf{x}_i \left(\tau\right), \tau ; \mathbf{X}\right)$ with $ \displaystyle \fint_{B_{\delta}(\mathbf{x}_i(\tau))}\nabla f \left(\xi, \tau ; \mathbf{X}\right) d\xi$.\\
We shall sketch the proof of the result looking back to the proof of Theorem 1, analyzing only the main points affected by such a change.\\
The mapping defined in (\ref{operatore}) for (\ref{dinamica}) can be rewritten likewise for (\ref{dinamica_media}), 
\begin{equation}\label{operatore2}
\Psi\left(\mathbf{X}, \mathbf{V}\right) (t) = 
\begin{pmatrix}
{\mathbf{x}_0}_1+ \displaystyle \int_{0}^{t}\mathbf{v}_1 (\tau) \, d \tau \\
\vdots \\
{\mathbf{x}_0}_n+ \displaystyle \int_{0}^{t}\mathbf{v}_n (\tau) \, d \tau \\
{\mathbf{v}_0}_1+ \displaystyle \int_{0}^{t}\mathbf{F}_1 \left( \tau, \mathbf{X}(\tau), \mathbf{V}(\tau),\fint_{B_{\delta}(\mathbf{x}_1(\tau))} \nabla f\left( \xi,\tau;\mathbf{X}\right) d\xi \right) \, d \tau \\
\vdots \\
{\mathbf{v}_0}_n+ \displaystyle \int_{0}^{t}\mathbf{F}_n \left( \tau, \mathbf{X}(\tau), \mathbf{V}(\tau), \fint_{B_{\delta}(\mathbf{x}_n(\tau))} \nabla f\left( \xi,\tau;\mathbf{X}\right) d\xi\ \right) \, d \tau \\
\end{pmatrix}
\end{equation}
for any $ \left( \mathbf{X}, \mathbf{V}\right) \in E_{\widetilde{R}} $, $t \in [0,\widetilde{T}]$, where $E_{\widetilde{R}}:=C\left( [0,\widetilde{T}];  B_{\widetilde{R}}(\mathbf{X}_0) \times B_{\widetilde{R}}(\mathbf{V}_0) \right)$, for a fixed $\widetilde{R}>0$, $0<\widetilde{T}\le T$.\\
Let $ \left(\mathbf{X}, \mathbf{V} \right)$ and $\left( \widehat{\mathbf{X}}, \widehat{\mathbf{V}} \right)  \in E_{\widetilde{R}}$.
%
Clearly ($\ref{contra_prime}$) still remains true. We now focus on the component $l \in \{n+1,...,2n\}$.
 From the Lipschitz condition ($\ref{Flip}$), it follows that 
\begin{equation}\label{contra_second_media}
\left.
\begin{array}{rl}
&\left |  \Psi_{n+j}\left(\mathbf{X}, \mathbf{V}\right) (t)  - \Psi_{n+j}\left(\widehat{\mathbf{X}}, \widehat{\mathbf{V}}\right) (t) \right | \le\\\\
&\le L_{F}^{R} \displaystyle \int_{0}^{t} \left |  \mathbf{X}(\tau)- \widehat{\mathbf{X}}(\tau)\right | + \left |  \mathbf{V}(\tau)- \widehat{\mathbf{V}}(\tau)\right | d\tau +\\\\
\ \ \ \ &+ L_{F} \displaystyle \int_{0}^{t} \left | \displaystyle \fint_{B_{\delta}(\mathbf{x}_j(\tau))}\nabla f \left( \xi, \tau ; \mathbf{X}\right) d\xi  - \fint_{B_{\delta}(\mathbf{\hat{x}}_j(\tau))}\nabla f \left( \xi, \tau ; \widehat{\mathbf{X}}\right) d\xi \right| d\tau \le \\\\
& \le L_{F}^{R}  \overline{T} \left(  || \mathbf{X} - \widehat{\mathbf{X}}||_{\infty,\widetilde{T}}+ || \mathbf{V} - \widehat{\mathbf{V}}||_{\infty, \widetilde{T}}\right) + L_{F} \left( J_1 + J_2 \right),
\end{array}
\right.
\end{equation}
where $ J_1, J_2$ are the integrals
\begin{equation}
J_1 :=   \displaystyle \int_{0}^{t} \left | \displaystyle \fint_{B_{\delta}(\mathbf{x}_j(\tau))}\nabla f \left( \xi, \tau ; \mathbf{X}\right) d\xi  - \fint_{B_{\delta}(\mathbf{\hat{x}}_j(\tau))}\nabla f \left( \xi, \tau ; \mathbf{X}\right) d\xi \right| d\tau,
\end{equation}
\begin{equation}
J_2 :=\displaystyle \int_{0}^{t} \left | \displaystyle \fint_{B_{\delta}(\mathbf{\hat{x}}_j(\tau))}\nabla f \left( \xi, \tau ; \mathbf{X}\right) d\xi  - \fint_{B_{\delta}(\mathbf{\hat{x}}_j(\tau))}\nabla f \left( \xi, \tau ; \widehat{\mathbf{X}}\right) d\xi \right| d\tau.
\end{equation}
Since $\left| B_{\delta}(\mathbf{x}_j(\tau)) \right|=\left| B_{\delta}(\mathbf{\hat{x}}_j(\tau)) \right|$, we obtain:
\begin{equation}\label{J_1}
\left.
\begin{array}{rl}
&\left | \displaystyle \fint_{B_{\delta}(\mathbf{x}_j(t))}\nabla f \left( \xi, \tau ; \mathbf{X}\right) d\xi  - \fint_{B_{\delta}(\mathbf{\hat{x}}_j(t))}\nabla f \left( \xi, \tau ; \mathbf{X}\right) d\xi \right| =\\\\
& \le \displaystyle \frac{1}{\left| B_{\delta}(\mathbf{x}_j(\tau)) \right|}\left | \int_{B_{\delta}(\mathbf{x}_j(\tau))}\nabla f \left( \xi, \tau ; \mathbf{X}\right) d\xi -  \int_{B_{\delta}(\mathbf{\hat{x}}_j(\tau))}\nabla f \left( \xi, \tau ; \mathbf{X}\right) d\xi\right| \le \\\\
& \le \displaystyle \frac{1}{\delta^{N}C_N} \displaystyle \int_{B_{\delta}} \left| \nabla f \left( \xi - \mathbf{x}_j(\tau) , \tau ; \mathbf{X}\right)-  \nabla f \left( \xi - \mathbf{\hat{x}}_j(\tau) , \tau ; \mathbf{X}\right) \right| d\xi,
\end{array}
\right.
\end{equation}
where $C_N$ denotes the volume of $B_{\delta} \subset \mathbb{R}^N$.\\
The estimate follows as in ($\ref{stimI1}$), since from ($\ref{valor_medio}$) we get:
\begin{equation}\label{stimI1_media}
\left.
\begin{array}{lll}
J_1 &\le \displaystyle N^2 K e^{\left(\delta^{2}+|\mathbf{X}_{0}|^{2}+\widetilde{R}^{2}\right)} \int_{0}^{t} \left(\frac{H}{\tau^{1-\alpha/2}}+\frac{2 H_{\bf X}}{\alpha}\tau^{\alpha/2}\right) \left| \mathbf{x}_{j}(\tau)- \widehat{\mathbf{x}}_{j}(\tau)\right| d \tau \le \\\\
& \le N^2 K e^{\left(\delta^{2}+|\mathbf{X}_{0}|^{2}+\widetilde{R}^{2}\right)} \left( \displaystyle \frac{2H}{\alpha} \displaystyle\widetilde{T}^{\alpha/2} +\displaystyle \frac{2H_{\bf X}}{\alpha+2}\displaystyle \widetilde{T}^{\alpha/2+1}  \right) || \mathbf{X} - \widehat{\mathbf{X}}||_{\infty,\widetilde{T}}.
\end{array}
\right.
\end{equation}
Recalling the estimates of $I_2$ in the previous section, we observe that 
\begin{equation}\label{stimJ2_media}
\left.
\begin{array}{lll}
J_2 &=\displaystyle \int_{0}^{t} \left | \displaystyle \fint_{B_{\delta}(\mathbf{\hat{x}}_j(\tau))}\nabla f \left( \xi, \tau ; \mathbf{X}\right)  - \nabla f \left( \xi, \tau ; \widehat{\mathbf{X}}\right) d\xi \right| d\tau\le \\\\
& \le \displaystyle \int_{0}^{t} \frac{1}{\left| B_{\delta}(\mathbf{\hat{x}}_j(\tau))\right|} \int_{B_{\delta}(\mathbf{\hat{x}}_j(\tau))}\left|\nabla f \left( \xi, \tau ; \mathbf{X}\right)  - \nabla f \left( \xi, \tau ; \widehat{\mathbf{X}}\right)  \right| d\xi d\tau \le\\\\
& \le \displaystyle \int_{0}^{t} \frac{1}{\left| B_{\delta}(\mathbf{\hat{x}}_j(\tau))\right|} \left(\int_{B_{\delta}(\mathbf{\hat{x}}_j(\tau))}\left|\nabla \eta \left( \xi, \tau \right)   \right| d\xi \right)d\tau ,
\end{array}
\right.
\end{equation}
where $\eta$ has been introduced in ($\ref{peqdiff}$).
From ($\ref{stim3}$) we get
\begin{equation}\label{stimJ2_finale}
\left.
\begin{array}{rl}
J_{2} &\le \displaystyle C_{\Gamma} H_{\widetilde{R}+\delta} e^{2C\left(\left|  \mathbf{X}_0\right|^2+\widetilde{R}^2 +\delta^2\right)}  \displaystyle \left(\frac{\pi}{\overline{\gamma}}\right)^{N/2}  \int_{0}^{t} \sqrt{\tau} d\tau\left | \left |  \mathbf{X}- \widehat{\mathbf{X}} \right|\right |_{\infty,\widetilde{T}}\le\\\\
&\le  C_{\Gamma} H_{\widetilde{R}+\delta} e^{2C\left(\left|  \mathbf{X}_0\right|^2+\widetilde{R}^2+\delta^2\right)}   \displaystyle  \left(\frac{\pi}{\overline{\gamma}}\right)^{N/2} \widetilde{T}^{\frac{3}{2}}\left | \left |  \mathbf{X}- \widehat{\mathbf{X}} \right|\right |_{\infty,\widetilde{T}}.
\end{array}
\right.
\end{equation}
From (\ref{stimI1_media}) and (\ref{stimJ2_finale}) we deduce a sufficient condition, of the same form as (\ref{minT}), that $\widetilde{T}$ has to satisfy in order to ensure that $\Psi$ is a contraction operator.
As done in (\ref{secondecomp}), we can easily show an analogous relation for $\widetilde{T}$ in order to guarantee that the range of  $\Psi $ is a subset of $E_{\widetilde{R}}$. In fact, we observe that the replacement of the gradient term with the average results in a slight modification of the exponential factor of $(\ref{comp_prov})$, that is $e^{\kappa \left( |\mathbf{X}_0|^{2}+R^{2}\right)}$ is replaced by $e^{\kappa \left( |\mathbf{X}_0|^{2}+\widetilde{R}^{2}+ \delta^2 \right)}$.
%
\section{Existence of global solution result.}
\noindent
In this section we present a global existence and uniqueness result for (\ref{dinamica})-(\ref{peq}) using a principle of continuation of solutions. Strengthening the growth conditions on functions $g$ and $\varphi$, we shall prove that bounded solutions can be continued.
For the sake of completeness, an analogous result will be also shown for (\ref{dinamica_media})-(\ref{peq}). \\
In light of the previous section, let assume that there exists a local solution $\mathbf{Y}=\mathbf{Y}(t)$  to (\ref{dinamica})-(\ref{peq}) in $[0,\overline{T}]$, with $\overline{T} \le T$. 
We replace the previous hypotheses H1) and H4) 
 with the following assumptions:
\begin{itemize}
\item[H6)] For every $ i=1,...,n$
\begin{equation}\label{hp_F}
\mathbf{F}_i:  \left[0,T\right] \times \mathbb{R}^{N\times n} \times \mathbb{R}^{N\times n} \times  \mathbb{R}^{N} \rightarrow \mathbb{R}^{N}
\end{equation}
is globally Lipschitz continuous, with Lipschitz constant $L_F$.
\item[H7)]
$g: \mathbb{R}^N \times \mathbb{R}^{N \times n} \longrightarrow  \mathbb{R}$, $\varphi: \mathbb{R}^N  \longrightarrow  \mathbb{R}$
are assumed to be a continuous functions, satisfying:
\begin{equation}\label{cresc_lin2}
\left| g(x,\mathbf{X}) \right| \le M(1+\left|x \right|+  \left|\mathbf{X} \right|), \ \ \ \ \ \left| \varphi(x) \right| \le M(1+\left|x \right|)
\end{equation}
for any $ x \in \mathbb{R}^{N}$, $\forall \mathbf{X} \in C\left(\left[0,\overline{T}\right];\mathbf{R}^{N\times n}\right)$, with $\overline{T}\le T$.\\
\end{itemize}
\begin{Theorem}\label{teoglobale}\textbf{Global existence and uniqueness of the solution}\\
Under assumptions H6)-H7), the local solution $\mathbf{Y}$ is a global solution.
\end{Theorem}
In order to state Theorem \ref{teoglobale}, we shall use the following result:
\begin{Lemma}\label{lemmagrad}
Let $\mathbf{X} \in C\left(\left[0,\overline{T}\right];\mathbf{R}^{N\times n}\right)$, with $\overline{T}\le T$.\\
Under assumption H7), the following estimate holds true:
\begin{equation}\label{lemma_exp}
\left| \nabla f \left( x,t; \mathbf{X}\right) \right|\le K_1 \left( \frac{1+ \left| x\right|}{\sqrt{t}}+ K_2+ \int_0^{t} \left(\frac{1+|x|+ |\mathbf{X}(\tau)|}{\sqrt{t-\tau}}\right)d\tau \right),
\end{equation}
for any $ x \in \mathbf{R}^{N}$, $t \in \left(0, \overline{T} \right]$,
where $K_1=K_1(C_{\Gamma},N,M,\lambda_0^{*})$, $K_2=K_2(N,\lambda_0^{*},T)$, $\lambda_0^{*}<\lambda_0$, $C_{\Gamma}$ as in (\ref{der_gamma}). \\
\end{Lemma}
Moreover we adopt the following Gronwall-type inequality:
\begin{Lemma}\label{lemma2}
Let $h \in C \left( [0,T]\right)$, $w\in L^1 \left( [0,T]\right)$, $v \in L^1 \left( [0,T] \times [0,T]\right)$ be non-negative functions, such that 
\begin{equation}\label{formagron1}
h(t) \le \alpha + \displaystyle \int_{0}^{t} w(\tau)h(\tau) d\tau+  \displaystyle \int_{0}^{t} \int_{0}^{\tau} v(s,\tau)h(s) ds d\tau \ \ \ \ \forall t \in [0,T],
\end{equation}
Then the following inequality holds true:
\begin{equation}
h(t) \le   \displaystyle \alpha  \ \ \mbox{exp}{ \left[ \displaystyle \int_{0}^{t} \left( w(\tau)+ \displaystyle  \int_{0}^{\tau}v(s,\tau) ds \right) d\tau\right]} \ \ \ \ \forall t \in [0,T].
\end{equation}

\end{Lemma}
%
\noindent
\textit{Proof of Theorem \ref{teoglobale}.} To establish the result, we structure the proof in two steps. In the first we shall prove that a solution $\mathbf{Y}(t)$ of (\ref{dinamica})-(\ref{peq}) on an interval $[0,\overline{T})$ remains bounded. In the second, we will prove that $\mathbf{Y}(t)$ can be extended to $[0,T]$, showing that there exists the $\lim_{t \rightarrow\overline{T}^{-} }\mathbf{Y}(t)$. Thesis follows applying Theorem \ref{teolocale} starting at $t=\overline{T}$.
\textit{Step 1:}
Recalling $(\ref{dinamica})_{2}$ we get
$\forall i=1,...,n$
\begin{equation}
\mathbf{v}_i(t)={\mathbf{v}_0}_i+ \displaystyle \int_{0}^{t}\mathbf{F}_i \left( \tau, \mathbf{X}(\tau), \mathbf{V}(\tau), \nabla f\left( \mathbf{x}_i(\tau), \tau; \mathbf{X} \right)\right) \, d \tau ,\\
\end{equation}
thus
\begin{equation}\label{stimVi}
\left.
\begin{array}{rl}
&\left|\mathbf{v}_i(t)-{\mathbf{v}_0}_i\right| \le \displaystyle \int_{0}^{t} \left| \mathbf{F}_i \left( \tau, \mathbf{X}(\tau), \mathbf{V}(\tau), \nabla f\left( \mathbf{x}_i(\tau), \tau; \mathbf{X} \right)\right) \right| \, d \tau \le\\\\
&\le t\displaystyle C_0 + L_F \displaystyle \int_{0}^{t} \left| \mathbf{X}(\tau)-  \mathbf{X}_0 \right| d\tau + \\\\
& + L_F \displaystyle \int_{0}^{t} \left| \mathbf{V}(\tau)-  \mathbf{V}_0 \right| d\tau+
 L_F \displaystyle \int_{0}^{t} \left| \nabla f \left( \mathbf{x}_i(\tau), \tau; \mathbf{X} \right) \right| d\tau
\end{array}
\right.
\end{equation}
From (\ref{lemma_exp}), we get 
\begin{equation}\label{stimV}
\left.
\begin{array}{rl}
&\left|\mathbf{V}(t)-{\mathbf{V}_0}\right| \le n T C_0 + nL_F \displaystyle \int_{0}^{t} \left| \mathbf{X}(\tau)-  \mathbf{X}_0 \right| d\tau + \\\\
& + nL_F \displaystyle \int_{0}^{t} \left| \mathbf{V}(\tau)-  \mathbf{V}_0 \right| d\tau+
 nL_F K_1\left( \int_{0}^{t}\frac{1+\left|\mathbf{X}(\tau)\right|}{\sqrt{\tau}} d\tau+K_2T\right) +\\\\
 &+nL_F K_1\displaystyle  \left( \int_{0}^{t}\int_{0}^{\tau}\frac{1+\left|\mathbf{X}(\tau)\right|+\left|\mathbf{X}(s)\right|}{\sqrt{\tau-s}} ds d\tau \right).
\end{array}
\right.
\end{equation}

Integrating $(\ref{dinamica})_{1}$, we get
\begin{equation}\label{stimaX}
\mathbf{X}(t)=\mathbf{X}_0+ \displaystyle \int_{0}^{t} \mathbf{V}(\tau) \, d \tau = \mathbf{X}_0+\mathbf{V}_0 t+ \displaystyle \int_{0}^{t} \left| \mathbf{V}(\tau)-  \mathbf{V}_0 \right| d\tau.
\end{equation}
By equations (\ref{stimaX}) and (\ref{stimV}) we obtain that
\begin{equation}\label{stima_completa}
\left.
\begin{array}{rl}
&\left|\mathbf{Y}(t)-\mathbf{Y}_0\right| \le \left|\mathbf{X}(t)-\mathbf{X}_0\right|+  \left|\mathbf{V}(t)-\mathbf{V}_0\right| \le \left| \mathbf{V}_0 \right| T+ \displaystyle \int_{0}^{t} \left| \mathbf{Y}(\tau)-  \mathbf{Y}_0 \right| d\tau+\\\
&+n T C_0 + nL_F\sqrt{2} \displaystyle \int_{0}^{t} \left| \mathbf{Y}(\tau)-  \mathbf{Y}_0 \right| d\tau + \\\\
&+ nL_FK_1\left(  \left(1+  \left|\mathbf{X}_0\right|\right)2\sqrt{T}+ \displaystyle \int_{0}^{t} \frac{\left| \mathbf{Y}(\tau)-  \mathbf{Y}_0 \right|}{\sqrt{\tau}} d\tau + K_2 T \right)+\\\\
&+ nL_FK_1 \left( \left(1+  2\left|\mathbf{X}_0\right|\right)\frac{4}{3}\sqrt{T^3}+ \displaystyle \int_{0}^{t} \int_{0}^{\tau} \frac{\left| \mathbf{X}(\tau)-  \mathbf{X}_0 \right|}{\sqrt{\tau-s}} ds d\tau + \right.\\\\
&+ \left. \displaystyle \int_{0}^{t} \int_{0}^{\tau} \frac{\left| \mathbf{X}(s)-  \mathbf{X}_0 \right|}{\sqrt{\tau-s}} ds d\tau \right)\le \\\\
&\le \left| \mathbf{V}_0 \right| T +n T C_0+nL_FK_1\left(1+  \left|\mathbf{X}_0\right|\right)2\sqrt{T}+nL_FK_1\left(1+  2\left|\mathbf{X}_0\right|\right)\frac{4}{3}\sqrt{T^3}+\\\\
&+nL_FK_1 K_2T + \left(1+nL_F\sqrt{2}\right) \displaystyle \int_{0}^{t} \left| \mathbf{Y}(\tau)-  \mathbf{Y}_0 \right| d\tau+\\\\
&+nL_FK_1  \displaystyle \int_{0}^{t} \frac{\left| \mathbf{Y}(\tau)-  \mathbf{Y}_0 \right|}{\sqrt{\tau}} d\tau+
nL_FK_1  \displaystyle \int_{0}^{t}\int_{0}^{\tau} \frac{\left| \mathbf{Y}(\tau)-  \mathbf{Y}_0 \right|}{\sqrt{\tau-s}}ds d\tau+\\\\
&+nL_FK_1  \displaystyle \int_{0}^{t}\int_{0}^{\tau} \frac{\left| \mathbf{Y}(s)-  \mathbf{Y}_0 \right|}{\sqrt{\tau-s}}ds d\tau.
 \end{array}
\right.
\end{equation}
Denoting with $\alpha=\alpha\left( \mathbf{X}_0, \mathbf{V}_0, L_F, K_1, K_2, N,n,T\right)$ the quantity $\left| \mathbf{V}_0 \right| T +n T C_0+nL_FK_1\left(1+  \left|\mathbf{X}_0\right|\right)2\sqrt{T}+nL_FK_1\left(1+  2\left|\mathbf{X}_0\right|\right)\frac{4}{3}\sqrt{T^3}+nL_FK_1 K_2T $, we rewrite
\begin{equation}\label{pergron}
\left.
\begin{array}{rl}
&\left|\mathbf{Y}(t)-\mathbf{Y}_0\right| \le 
\alpha+ \left(1+nL_F\sqrt{2}\right) \displaystyle \int_{0}^{t} \left| \mathbf{Y}(\tau)-  \mathbf{Y}_0 \right| d\tau+\\\\
&+nL_FK_1  \displaystyle \int_{0}^{t} \frac{\left| \mathbf{Y}(\tau)-  \mathbf{Y}_0 \right|}{\sqrt{\tau}} d\tau+
nL_FK_1  \displaystyle \int_{0}^{t}\left| \mathbf{Y}(\tau)-  \mathbf{Y}_0 \right|2\sqrt{\tau} d\tau+\\\\
&+nL_FK_1  \displaystyle \int_{0}^{t}\int_{0}^{\tau} \frac{\left| \mathbf{Y}(s)-  \mathbf{Y}_0 \right|}{\sqrt{\tau-s}}ds d\tau=\\\\
&= \alpha+  \left(1+nL_F\sqrt{2}\right) \displaystyle \int_{0}^{t} \left| \mathbf{Y}(\tau)-  \mathbf{Y}_0 \right| d\tau+nL_FK_1  \displaystyle \int_{0}^{t} \frac{\left| \mathbf{Y}(\tau)-  \mathbf{Y}_0 \right|}{\sqrt{\tau}} d\tau+\\\\
&+2nL_FK_1 \displaystyle \int_{0}^{t}\sqrt{\tau}\left| \mathbf{Y}(\tau)-  \mathbf{Y}_0 \right|d\tau+nL_FK_1  \displaystyle \int_{0}^{t}\int_{0}^{\tau} \frac{\left| \mathbf{Y}(s)-  \mathbf{Y}_0 \right|}{\sqrt{\tau-s}}ds d\tau.
 \end{array}
\right.
\end{equation}
Finally, we observe that (\ref{pergron}) can be rewritten as

\begin{equation}\label{formagron}
h(t) \le \alpha + \displaystyle \int_{0}^{t} w(\tau)h(\tau) d\tau+  \displaystyle \int_{0}^{t} \int_{0}^{\tau} v(s,\tau)h(s) ds d\tau,
\end{equation}
with $h(t):=\left|\mathbf{Y}(t)-\mathbf{Y}_0\right| \in C\left(0,T\right)$, and $w,v \in L^1\left([0,T]\right)$ defined by $w(t):=1+nL_F\sqrt{2}+ \displaystyle \frac{nL_FK_1}{\sqrt{t}}+2nL_FK_1\sqrt{t} $,\\
\begin{equation}
 v(t,\tau):=
\left\{
\begin{array}{rl}
\displaystyle \frac{nL_FK_1}{\sqrt{t-\tau}} \ \ \ \ \tau < t  \\\
 0 \ \ \ \ \ \ \ \ \tau \ge t.
 \end{array}
\right.
\end{equation}
\medskip
\medskip
Thus, by Lemma \ref{lemma2}, we reach the conclusion 
\begin{equation}\label{bound}
\left.
\begin{array}{rl}
\left|\mathbf{Y}(t)-\mathbf{Y}_0\right| \le  \alpha \ \ \mbox{exp} \left[ \displaystyle \int_0^t \left ( w(\tau)+ \int_0^{\tau} v(s,\tau)ds \right) d\tau\right]:= B 
 \end{array}
\right.
\end{equation}
where we denote with $B$ the following constant:
\begin{equation}
\left.
\begin{array}{rl}
B &= B\left( n,\alpha, L_F, M,C_{\Gamma},\lambda_{0}^{*},N, T\right):=\\
&=\alpha \ \ \mbox{exp} \left[\left( 1+nL_F \sqrt{2}\right)T+2\sqrt{T}nL_FK_1+ \displaystyle \frac{2}{3}nL_FK_1T^{3/2}\right]=\\
&= \alpha \ \ \mbox{exp} \left[\left( 1+nL_F \sqrt{2}\right)T+2\sqrt{T}nL_F C_{\Gamma} M  \displaystyle \frac{2^N \pi^{\frac{N}{2}}}{\left( \lambda_0^{*}\right)^{\frac{N}{2}}} \left( 1+\frac{T}{3}\right) \right] . 
 \end{array}
\right.
\end{equation}
%
\textit{Step 2:}
In order to prove the existence of a limit for $\left|\mathbf{Y}(t)\right|$, as $t$ goes to $T^{-}$, we consider $\{t_n\}$ a monotonic increasing sequence, with $\lim_{t\rightarrow T}$$\{ t_n\}$, and show that  $\mathbf{Y}(t_n)$ is a Cauchy sequence.\\
Let $0<t_m<t<t_n<T$. 
From (\ref{bound}) follows, in particular, that $\left|\mathbf{X}(t)\right| \le B+\left| \mathbf{X}_0\right|$, $\left|\mathbf{V}(t)\right| \le B + \left| \mathbf{V}_0\right|$, for any $t \in [0,\overline{T}]$. 
This allow us to retrace the computation of Step 1, estimating $\left|\mathbf{X}_0\right|$, $\left|\mathbf{V}_0\right|$ with $2B$, and obtaining that
\begin{equation}
\left.
\begin{array}{rl}
\displaystyle \left|\mathbf{Y}(t)-\mathbf{Y}(t_m)\right| \le A(t-t_m)+ \int_{t_m}^{t}w(\tau)\left|\mathbf{Y}(\tau)-\mathbf{Y}(t_m)\right| d\tau+\\\\
\displaystyle \int_{t_m}^{t} \int_{t_m}^{\tau}v(s,\tau)\left|\mathbf{Y}(s)-\mathbf{Y}(t_m)\right| ds d\tau,
\end{array}
\right.
\end{equation}
where 
$
A(t)=A\left(B, L_F, N,n,t\right):= \displaystyle 2B t +n\left(C_0+L_FK_1K_2 \right) t+nL_FK_1\left(1+2B\right)2\sqrt{t}+nL_FK_1\left(1+ 4B\right)\frac{4}{3}\sqrt{t^3}.
$\\
By Lemma \ref{lemma2} and (\ref{bound}) we get that
\begin{equation}
\left|\mathbf{Y}(t_n)-\mathbf{Y}(t_m)\right| \le A(t_n-t_m) \mbox{exp}\left[  \int_{t_m}^{t_n} \left( w(\tau)+ \int_{t_m}^{\tau}v(s,\tau)ds \right) d\tau\right] 
= A(t_n-t_m) \ \displaystyle \frac{B}{\alpha}.
\end{equation}
From the expression of $A$ we observe that $A(t_n-t_m) \rightarrow 0$ as $n,m \rightarrow \infty$.\\
Thus $\left \{\mathbf{Y}(t_n)\right \}$ is a Cauchy sequence, and admits a limit valute as $t \rightarrow T^{-}$.\\

\noindent
\textit{Proof of Lemma  \ref{lemmagrad}.}
From (\ref{fgamma}) we get the following expression
\begin{equation}
\left.
\begin{array}{rl}
\nabla f\left(x,t; \mathbf{X}\right)&= \displaystyle \int_{\mathbb{R}^{N}} \nabla \Gamma \left( x,t;\xi,0\right) \varphi(\xi) d\xi - \int_{0}^{t}\int_{\mathbb{R}^{N}} \nabla \Gamma \left( x,t;\xi,\tau\right)g_{\bf X}(\xi,\tau) d\xi d\tau = \\\\
&= \widetilde{G_1}-\widetilde{G_2},
\end{array}
\right.
\end{equation}
where
\begin{equation}
\widetilde{G_1}= \displaystyle \int_{\mathbb{R}^{N}} \nabla \Gamma \left( x,t;\xi,0\right) \varphi(\xi) d\xi,\ \ \ \
\widetilde{G_2}=  \displaystyle \int_{0}^{t}\int_{\mathbb{R}^{N}} \nabla \Gamma \left( x,t;\xi,\tau\right)g_{\bf X}(\xi,\tau) d\xi d\tau.
\end{equation}
We now focus on $\widetilde{G_1}$. In particular, from (\ref{der_gamma}) and (\ref{cresc_lin2}) we obtain
\begin{equation}\label{stimlin}
\left.
\begin{array}{rl}
&\left|  \displaystyle \widetilde{G_1} \right| \le  \displaystyle \int_{\mathbb{R}^{N}} \displaystyle C_{\Gamma} \frac{1}{t^{\frac{N+1}{2}}} \displaystyle e^{ -\frac{\lambda_0^{*}}{4}\frac{\left|\xi-x\right|^2}{t}} M(1+\left| \xi \right|) d\xi \le \\\\
&\le \displaystyle C_{\Gamma} M  \frac{1}{\sqrt{t}} \displaystyle \int_{\mathbb{R}^{N}} 
 e^{ -\frac{\lambda_0^{*}}{4} \left| u \right|^{2}}  du \left( 1+ \left| x \right|\right)+ C_{\Gamma} M \displaystyle \int_{\mathbb{R}^{N}} 
 e^{ -\frac{\lambda_0^{*}}{4} \left| u \right|^{2}} \left| u \right| du=\\\\
 &=\displaystyle C_{\Gamma} M  \frac{1}{\sqrt{t}} \displaystyle I_{0}\left(\frac{\lambda_0^*}{4}\right) \left( 1+ \left| x \right|\right)+ C_{\Gamma} M \displaystyle I_{1}\left(\frac{\lambda_0^*}{4}\right),
\end{array}
\right.
\end{equation}
where, in the last inequality, $I_0(.)$, $I_1(.)$ are defined in $(\ref{I0app})$, $(\ref{I1app})$.\\
Replacing $\displaystyle I_{0}\left(\frac{\lambda_0^*}{4}\right)$ and $\displaystyle I_{1}\left(\frac{\lambda_0^*}{4}\right)$ with their exact value (see Appendix for computation), we obtain the following estimate for $\widetilde{G_1}$ :
\begin{equation}\label{stimlin2}
\left.
\begin{array}{rl}
&\left|  \displaystyle \widetilde{G_1} \right| \le   \displaystyle C_{\Gamma} M \left( \frac{2^N \pi^{\frac{N}{2}}}{\left( \lambda_0^{*}\right)^{\frac{N}{2}}} \frac{\left( 1+ \left| x \right|\right)}{\sqrt{t}}+ \left(\frac{4}{\lambda_0^{*}}\right)^{\frac{N+1}{2}}\frac{\omega_N}{2}\frac{2\pi^{\frac{N+1}{2}}}{\omega_{N+1}}\right)=\\\\
&= \displaystyle C_{\Gamma} M \left( \frac{2^N \pi^{\frac{N}{2}}}{\left( \lambda_0^{*}\right)^{\frac{N}{2}}} \frac{\left( 1+ \left| x \right|\right)}{\sqrt{t}}+ \frac{2^{N}2}{\left(\lambda_0^{*}\right)^{\frac{N}{2}}\sqrt{\lambda_0^{*}}}\frac{\omega_N\pi^{\frac{N}{2}}\pi}{\omega_{N+1}}\right)=\\\\
&= \displaystyle C_{\Gamma} M  \frac{2^N \pi^{\frac{N}{2}}}{\left( \lambda_0^{*}\right)^{\frac{N}{2}}}
\left(\frac{\left( 1+ \left| x \right|\right)}{\sqrt{t}}+ \frac{2}{\sqrt{\lambda_0^{*}}}\frac{\omega_N\pi}{\omega_{N+1}}\right)=K_1\left(\frac{\left( 1+ \left| x \right|\right)}{\sqrt{t}}+\widetilde{K}_2\right),
\end{array}
\right.
\end{equation}
where $K_1:=C_{\Gamma} M  \displaystyle \frac{2^N \pi^{\frac{N}{2}}}{\left( \lambda_0^{*}\right)^{\frac{N}{2}}}$, and $\widetilde{K}_2:= \displaystyle \frac{2}{\sqrt{\lambda_0^{*}}}\frac{\omega_N\pi}{\omega_{N+1}}$.
As already noted in the proof of Proposition 1, the estimate of $\widetilde{G_2}$ follows by similar computation. For the sake of simplicity, we omit details and state that
\begin{equation}\label{stimlin3}
\left.
\begin{array}{rl}
&\left|  \displaystyle \widetilde{G_2} \right| \le   \displaystyle K_1 \left( \int_{0}^{t}\frac{1+ \left| x \right| + \left| \mathbf{X}(\tau) \right|}{\sqrt{t-\tau}}d\tau+\widetilde{K}_2 T\right).
\end{array}
\right.
\end{equation}
Thesis follows from (\ref{stimlin2}) and (\ref{stimlin3}), denoting with $K_2$ the constant $\widetilde{K}_2 +\widetilde{K}_2 T$.


\noindent
\textit{Proof of Lemma \ref{lemma2}.}
Let us define $\widetilde{h}(t):= \displaystyle \sup_{0 \le s \le t} h(s)$. Clearly, $\widetilde{h}(t) \ge h(t)$, for any $ t \in [0,T]$. For any fixed $\widehat{t} \le T$, inequality (\ref{formagron}) implies that
\begin{equation}
h(t) \le \alpha + \displaystyle \int_{0}^{\widehat{t}} w(\tau)\widetilde{h}(\tau) d\tau+  \displaystyle \int_{0}^{t}\widetilde{h}(\tau) \int_{0}^{\tau} v(s,\tau) ds d\tau
\end{equation}
for any $t \le \widehat{t}$. Hence, in particular,
\begin{equation}
\widetilde{h}\left(\widehat{t}\right) \le \alpha + \displaystyle \int_{0}^{\widehat{t}} w(\tau)\widetilde{h}(\tau) d\tau+  \displaystyle \int_{0}^{\widehat{t}}\widetilde{h}(\tau) \int_{0}^{\tau} v(s,\tau) ds d\tau.
\end{equation}
We conclude that 
\begin{equation}\label{perlemma}
\widetilde{h}\left(t\right) \le \alpha + \displaystyle \int_{0}^{t} \left( w(\tau)+  \displaystyle \int_{0}^{\tau}v(s,\tau) ds \right) \widetilde{h}(\tau) d\tau \ \ \ \ \forall t \le T.
\end{equation}
Applying Gronwall inequality to (\ref{perlemma}), 
 we obtain
\begin{equation}
h(t)\le \widetilde{h}(t) \le \displaystyle \alpha  \ \ \mbox{exp}{ \left[ \displaystyle \int_{0}^{t} \left( w(\tau)+ \displaystyle  \int_{0}^{\tau}v(s,\tau) ds \right) d\tau\right]} \ \ \ \ \forall t \le T.
\end{equation}
This concludes the proof.\\
%

\noindent
\textbf{Remark 3.}
Under the same assumptions for the global existence of solutions to system (\ref{dinamica})-(\ref{peq}), we state an analogous result for system (\ref{dinamica_media})-(\ref{peq}), previously introduced as a modification of (\ref{dinamica}) and already investigated for a local result in Theorem 3.
\begin{Theorem}
Let $\mathbf{Y}=\mathbf{Y}(t)$ local solution to (\ref{dinamica_media})-(\ref{peq}) in $[0,\widetilde{T}]$, with $\widetilde{T} \le T$. \\
Under assumptions H6)-H7),  $\mathbf{Y}$ is a global solution.
 \end{Theorem}
\noindent
\textit{Proof.} We present only a sketch of the proof, which follows the same line of reasoning of proof of Theorem \ref{teoglobale}. In particular, we point out the steps in which the replacement of the gradient term in (\ref{dinamica}) with the average over a ball involves modifications. 
Integrating $(\ref{dinamica_media})_{2}$, we obtain that
\begin{equation}\label{stimVi_media}
\left.
\begin{array}{rl}
&\left|\mathbf{v}_i(t)-{\mathbf{v}_0}_i\right| \le \displaystyle \int_{0}^{t} \left|\mathbf{F}_i \left( \tau, \mathbf{X}(\tau), \mathbf{V}(\tau),\fint_{B_{\delta}(\mathbf{x}_i(t))} \nabla f\left( \xi,t;\mathbf{X}\right) d\xi\right) \right| \, d \tau \le\\\
&\le t\displaystyle C_0 + L_F \displaystyle \int_{0}^{t} \left| \mathbf{X}(\tau)-  \mathbf{X}_0 \right| d\tau + \\\\
& + L_F \displaystyle \int_{0}^{t} \left| \mathbf{V}(\tau)-  \mathbf{V}_0 \right| d\tau+
 L_F \displaystyle \int_{0}^{t} \left| \fint_{B_{\delta}(\mathbf{x}_i(t))} \nabla f\left( \xi,t;\mathbf{X}\right) d\xi \right| d\tau,
\end{array}
\right.
\end{equation}
The above Lemma \ref{lemmagrad}, recalling that $\xi \in B_{\delta}(\mathbf{x}_i(t)) $, leads to 
\begin{equation}\label{lemma_exp_media}
\left| \nabla f \left( \xi,\tau; \mathbf{X}\right) \right|\le K_1 \left( \frac{1+ \delta+|\mathbf{X}(\tau)|}{\sqrt{\tau}}+ K_2+ \int_0^{\tau} \left(\frac{1+\delta+  |\mathbf{X}(s)|+ |\mathbf{X}(\tau)|}{\sqrt{\tau-s}}\right)ds \right).
\end{equation}
Hence we get 
\begin{equation}\label{stimV_media}
\left.
\begin{array}{rl}
&\left|\mathbf{V}(t)-{\mathbf{V}_0}\right| \le n T C_0 + nL_F \displaystyle \int_{0}^{t} \left| \mathbf{X}(\tau)-  \mathbf{X}_0 \right| d\tau + \\\\
& + nL_F \displaystyle \int_{0}^{t} \left| \mathbf{V}(\tau)-  \mathbf{V}_0 \right| d\tau+
 nL_F K_1\left( \int_{0}^{t}\frac{1+\delta +\left|\mathbf{X}(\tau)\right|}{\sqrt{\tau}} d\tau+K_2T\right) +\\\\
 &+nL_F K_1\displaystyle  \left( \int_{0}^{t}\int_{0}^{\tau}\frac{1+\delta + \left|\mathbf{X}(\tau)\right|+\left|\mathbf{X}(s)\right|}{\sqrt{\tau-s}} ds d\tau\right).
\end{array}
\right.
\end{equation}
We note that the previous upper-bound for $\left|\mathbf{V}(t)-{\mathbf{V}_0}\right| $ can be obtained from (\ref{stimV}), simply replacing $\left|\mathbf{X}(\tau)\right|$ with $\left|\mathbf{X}(\tau)\right|+\delta$. Hence the proof of Theorem \ref{teoglobale} can be repeated for system (\ref{dinamica_media}), observing that $1+ \left| \mathbf{X_0} \right|$ will be replaced by $1+\delta+ \left| \mathbf{X_0} \right|$, and 
$1+2 \left| \mathbf{X_0} \right|$ by $1+\delta+ 2\left| \mathbf{X_0} \right|$, in (\ref{stima_completa}) and thus in the definition of  constant $\alpha$. We observe that modifications occur only in Step 1 of the proof, whereas Step 2 follows as for Theorem \ref{teoglobale}. 
%
\section{Conclusions and future works}
In this paper we present existence and uniqueness results of solutions for coupled hybrid systems of differential equations. 
From a modeling point of view, those kind of systems combine the advantages of individual-based models with continuous ones.  
The literature concerning analytical foundations is still lacking. Hence we stress the importance of our work, which has to be considered as a first step towards a more detailed analytical characterization of systems which are increasingly used to model phenomena.\\
One of the future work directions we are going to concerns another critical issue, that is the shortcoming of a detailed technique to estimate the parameters occurring in the models. 
Suppose to introduce a set of real parameters in a suitable domain $D \in \mathbb{R}^{p}$,  denoted with $\mathbf{\theta}=(\theta_1,...,\theta_p)$, in the reaction-diffusion equation. 
It means to assume $a_{ij}=a_{ij}( . ;\theta)$, $b_{i}=b_{i}( . ;\theta)$, $c=c( . ;\theta)$ in equation (\ref{intro3}), and $f=f( . ;\theta)$ as the solution to (\ref{intro3}), once assigned initial conditions.
Recalling that $f$ influences the dynamics expressed in (\ref{intro1}), the introduction of the parameters  affects the solution $\mathbf{X}=\mathbf{X}(\cdot;\theta)$. 
The introduction of parameters leads to the following modifications of equations (\ref{intro1})-(\ref{intro3}):
\begin{equation}\label{intro1para}
\mathbf{\ddot{x}}_i(t; \theta) = \displaystyle \mathbf{F}_{i}\left( t, \mathbf{X}(t; \theta),  \mathbf{\dot{X}}(t;\theta), \nabla f \left( \mathbf{x}_i(t;\theta),t;  \mathbf{X}(t;\theta), \theta \right),\theta\right),
\end{equation}
\begin{equation}\label{intro2para}
\left\{
\begin{array}{rl}
L^{\theta}f(x,t;  \mathbf{X}(t;\theta) )&=g(x,\mathbf{X}(t;\theta)) \ \ \ \ \ \ \ \ (x,t) \in \mathbb{R}^{N} \times (0,T),\\
f(x,0)&= \varphi(x) \ \ \ \ \ \ \ \ \ \ \ \ \ \ \ \ \ \ \ \ x \in \mathbb{R}^{N},
\end{array}
\right.
\end{equation}
where
\begin{equation}\label{intro3para}
L^{\theta} = \displaystyle \sum_{i,j=1}^{N} a_{i,j}(x,t;\theta) \partial_{x_i,x_j}^{2}+ 
	 \sum_{i=1}^{N} b_{i}(x,t;\theta) \partial_{x_i}+c(x,t;\theta)- \partial_{t}.
\end{equation}
Given a set of observable data,  \{$\mathbf{x}_{i}^{*}(t)$: $i=1,...,n$, $t \in I$ \} where $I$ is a discrete set of times, $I \subset (0,T)$, (e.g. the set of positions of $n$ cells recorded at different time steps of an experiment), the issue to address is to find the optimal values of the parameters for which a considered model fits reality, meaning that the model is able to reproduce the observed behavior.\\
Mathematically, this can be regarded as a least squares problem, aiming at finding $\theta_{opt} \in D$ arising from
$
\displaystyle \inf_{\theta \in D} \displaystyle \sum_{i=1}^{n}  \sum_{t \in I} \left| \mathbf{x}_{i}(.,\theta)- \mathbf{x}_{i}^{*}(t)\right|^{2}+ \epsilon(\theta),
$
where $ \epsilon(\theta)$ is a suitable penalization term, and $ \mathbf{x}_{i}(\cdot; \theta)$ is the solution to model (\ref{intro1para})-(\ref{intro2para}).
Standard numerical procedures used to solve the above optimization problem require the computation of the derivative ${\mathbf{u}_{i}}_{k}:=\partial_{\theta_k} \mathbf{x}_i \left(\cdot;\theta \right)$ for $k=1,...,p$.
By differentiating $(\ref{intro1para})$ and $(\ref{intro2para})$ with respect to $\theta_k$, $k=1,...,p$, we obtain that ${\mathbf{u}_{i}}_{k}$ is solution, at least formally, to a problem of the same form of (\ref{intro1para})-(\ref{intro2para}).

For that reason, we argue that analytical results presented in our paper can be extended to investigate the dependence of the solution on the parameters.
In a second ongoing work we investigate the case in which the spatial domain for (\ref{dinamica})-(\ref{peq}) is a bounded subset of $\mathbb{R}^N$. That will imply the introduction of suitable boundary conditions for (\ref{peq}).
Finally, we also extend our results to the case in which the source term of the parabolic equation is a less regular function, e.g. discontinuous.\\

\noindent
\textbf{Acknowledgements.} The authors extend warm thanks to Roberto Natalini, Director of Istituto per le applicazioni del calcolo ``Mauro Picone'', Consiglio Nazionale delle Ricerche (IAC-CNR), for many discussions and suggestions.

\section{Appendix}

\subsection{\textit{Fundamental solution and Cauchy problem}}
Let consider the Cauchy problem

\begin{equation}\label{cp_app}
\left\{
\begin{array}{rl}
Lu(x,t)&=f(x,t) \ \ \ \ \ \ \ \ (x,t) \in \Omega :=\mathbb{R}^N \times (0,T],\\
u(x,0)&= \varphi(x), \ \ \ \ \ \ \ \ \ \ \ x \in \mathbb{R}^N,
\end{array}
\right.
\end{equation}
where $L$ is the operator 
\begin{equation}\label{Lparabolic_app}
L = \displaystyle \sum_{i,j=1}^{N} a_{i,j}(x,t) \partial_{x_i,x_j}^{2}+ 
	 \sum_{i=1}^{N} b_{i}(x,t) \partial_{x_i}+c(x,t)- \partial_{t}.
\end{equation}
Let $L$ uniformly parabolic in $\Omega$, i.e. the matrix $\left( a_{ij}(x,t)\right)$ is positive definite and there exist positive constants $\overline{\lambda}_0$, $\overline{\lambda}_1$ such that for any $\xi \in \mathbb{R}^{n}$
$$
\overline{\lambda}_0\left| \xi \right|^{2} \le \sum_{i,j=1}^{n}  a_{i,j}(x,t)\xi_i\xi_j \le \overline{\lambda}_1\left| \xi \right|^{2} \ \ \ \ \ \ \ \ \ \forall (x,t) \in \Omega.
$$
Moreover let $a_{i,j}$, $b_{i}$, $c$ bounded H\"older continuous function in $\Omega$, with coefficient $\alpha \in (0,1)$ with respect to $x$ and $\alpha/2$ with respect to t.\\
Let $f(x,t)$, $\varphi(x)$ be continuous functions respectively in $\Omega$ and $\mathbb{R}^N$, satisfying
\begin{equation}
\left| f(x,t) \right| \le C e^{h |x|^2}, \  \ \ \ \left| \varphi(x) \right| \le C e^{h |x|^2},
\end{equation}
%
%
with $h$ positive constant satisfying $h < \displaystyle \frac{\lambda_0}{4T}$.\\
Finally, let assume $f$ locally H\"older continuous with exponent $\alpha$ in $x \in \mathbb{R}^N$, uniformly with respect to $t$.
Then the function
\begin{equation}
u(x,t) = \displaystyle \int_{ \mathbb{R}^N} \Gamma \left(x,t;\xi,0\right) \varphi\left(\xi\right)d \xi - \displaystyle  \int_{0}^{t}\int_{ \mathbb{R}^N} \Gamma \left(x,t;\xi,\tau\right)f\left(\xi, \tau\right) d \xi d\tau
\end{equation}
is a solution of (\ref{cp_app})-(\ref{Lparabolic_app}). \\
$\Gamma \left(x,t;\xi,\tau\right)$ is a fundamental solution of $Lu=0$,and it is a continuous function of $\left( x,t\right)$, uniformly with respect to $\left( \xi, \tau\right)$ if $t-\tau \ge constant>0$, and it is a continuous function of $\left( \xi,\tau \right)$, uniformly with respect to $\left( x,t\right)$ if $t-\tau \ge constant$. Hence $\Gamma \left(x,t;\xi,\tau\right)$ is a continuous function of $\left(x,t;\xi,\tau\right)$. Moreover also the first and second derivative with respect to space, and the derivative with respect to time of $\Gamma$ are continuous functions of $\left(x,t;\xi,\tau\right)$, where $x,\xi \in \mathbb{R}^{N}$ and $0 \le \tau < t \le T$ .
\subsection{\textit{Computation of a possible upper bound for $\lambda_0$ in Assumption H3).}}
Let A be a symmetric matrix satisfying
\begin{equation}\label{hip_A}
\mu_0 \left | \xi \right |^2 \le \left \langle A\xi, \xi \right \rangle
\le \mu_1 \left |\xi \right |^2 \ \ \ \ \forall \xi \in \mathbf{R}^{N}
\end{equation}
where $\mu_0,\mu_1$ are positive constants.\\
Defining $\eta=A\xi \in \mathbb{R}^{N}$,  $(\ref{unifpar})$ can be rewritten as

\begin{equation}\label{perlambda}
\mu_0 \left | A^{-1}\eta \right |^2 \le \left \langle \eta, A^{-1}\eta \right \rangle
\le \mu_1 \left |A^{-1}\eta \right |^2 .
\end{equation}
We choose $W=\{w_1,...,w_{N}\}$ an orthonormal basis of eigenvectors of A, and rewrite $\eta$ as a linear combination, with coefficients $c_1,...,c_N \in \mathbb{R}$,namely $\eta= \displaystyle \sum_{i=1}^{N} c_i w_i$.
Denoting with $\theta_i$, $i=1,...,N$, the eigenvalues of $A$, we thus obtain: 
$$
A^{-1}\eta= \displaystyle \sum_{i=1}^{N} \frac{c_i}{\theta_i} w_i.
$$\\
Since $\left | \eta \right |^{2}=\displaystyle \sum_{i=1}^{N} c_{i}^{2}$, we have

\begin{equation}
\left | A^{-1}\eta \right |^2=\left \langle A^{-1}\eta, A^{-1}\eta \right \rangle =  \displaystyle \sum_{i=1}^{N} \frac{c_{i}^{2} }{\theta_i^{2}} >  \frac{1}{ \displaystyle\max_{i=1,...,N} (\theta_i^{2})}  \displaystyle \sum_{i=1}^{N} c_{i}^{2}> \frac{1}{ \displaystyle \mu_1^{2}}  \left | \eta \right |^{2}.
\end{equation}
By $(\ref{perlambda})$, we get the following lower bound:

\begin{equation}\label{perlambda0}
\left \langle \eta, A^{-1}\eta \right \rangle \ge \mu_0 \left | A^{-1}\eta \right |^2 \ge \displaystyle \frac{\mu_0}{\mu_1^{2}}\left | \eta \right |^{2}.
\end{equation} 
We performed the algebraic computation above in order to show the relation between $\lambda_0$ and $\mu_0,\mu_1$ in H3), 
following the approach in \cite{friedman}.
Let $A$ be the matrix having entries $A_{ij}=a_{ij}(x,t)$, where $a_{ij}$ are the coefficients introduced in $(\ref{Lparabolic})$: we note that, since $L$ is of parabolic-type, $A$ is symmetric, and $(\ref{unifpar})$ corresponds to condition $(\ref{hip_A})$. Hence, from equation ($\ref{perlambda0}$), we get an upper bound for the parameter $\lambda_0$. 



\subsection{\textit{Computation of $I_0(\gamma)$, $I_1(\gamma)$}.}

Let $\gamma$ be a positive constant. In the following we give details of the computation of two integrals which appear in various points of the paper.

Let define $I_0\left( \gamma\right)$, $I_1\left( \gamma\right)$ the integrals
\begin{equation}\label{I0app}
I_0\left( \gamma\right):=\int_{\mathbb{R}^{N}} e^{ -\gamma \left| y \right|^{2}}  dy,
\end{equation}

\begin{equation}\label{I1app}
I_1\left( \gamma\right):=\int_{\mathbb{R}^{N}} e^{ -\gamma \left| y \right|^{2}} \left| y \right| dy.
\end{equation}

With the changes of variable given by $v:=\sqrt{2\gamma}y$ and then $u:=\frac{v}{\sqrt{2}}$, we rewrite:

\begin{equation}
I_0\left( \gamma\right):=\int_{\mathbb{R}^{N}} e^{ -\gamma \left| y \right|^{2}}  dy= \int_{\mathbb{R}^{N}} e^{ -\frac{v^2}{2}} \frac{1}{(2\gamma)^{N/2}} dv = \frac{1}{\gamma^{N/2}} \int_{\mathbb{R}^{N}} e^{ -  u^{2}}  du=\left(\frac{\pi}{\gamma}\right)^{N/2}.
\end{equation}

For the second, we convert it in polar coordinates and perform the change of variable given by $t=\gamma h^2$, obtaining 

\begin{equation}\label{stimlin}
\left.
\begin{array}{rl}
&I_1\left( \gamma\right):=\displaystyle \int_{\mathbb{R}^{N}} e^{ -\gamma \left| y \right|^{2}} \left| y \right|  dy
=\int_{0}^{\infty}e^{ -\gamma h^{2}}h^N \omega_N dh\\
&=\displaystyle{\int_{0}^{\infty}}e^{-t}\left(\frac{1}{\gamma}\right)^{\frac{N-1}{2}}t^{\frac{N-1}{2}}\frac{\omega_N}{2}\frac{1}{\gamma}dt=\displaystyle \Gamma_e\left(\frac{N+1}{2}\right)\left(\frac{1}{\gamma}\right)^{\frac{N+1}{2}}\frac{\omega_N}{2}=\\
&= \left(\frac{1}{\gamma}\right)^{\frac{N+1}{2}}\frac{\omega_N}{2}\frac{2\pi^{\frac{N+1}{2}}}{\omega_{N+1}},
\end{array}
\right.
\end{equation}
recalling that $\Gamma_e (N/2)\omega_N=2\pi^{\frac{N}{2}}$.\\

\newpage

\newpage

\section{Glossary}
\noindent
 Throughout the paper, we will identify each element of $\mathbb{R}^{N \times n}$ with the column vector of $\mathbb{R}^{Nn}$ obtained putting in column the $n$ columns one after another.\\

\begin{table}[h]
\label{tab:1}       
\begin{tabular}{lll}
$B_R(\mathbf{P})$ & Closed ball of $\mathbb{R}^{p} $, with center $\mathbf{P} \in \mathbb{R}^{p} $ and radius $R$ \\
$B_R$ & Closed ball of $\mathbb{R}^{p} $, with center $\mathbf{O} \in \mathbb{R}^{p} $ and radius $R$ \\
$|| .||_{\infty,T}$ &$\sup_{[0,T]}|.|$\\
$\omega_N$ & Surface area of the $(N-1)$-dimensional sphere, in $\mathbb{R}^N$\\ 
$\left| \mathbf{V}\right|$ & Euclidean norm of $\mathbf{V}$ $ \in \mathbb{R}^{p}$ \\
$\Gamma_e$ & Euler gamma function\\
$\mathbb{R}^{N \times n}$  & Space of matrices with $N$ rows and $n$ columns\\
$A=[\mathbf{a}_1...\mathbf{a}_d]$ & $A \in \mathbb{R}^{n \times d}$, with $\mathbf{a}_i \in \mathbb{R}^{n}$ $i=1,...,n$ column vector\\
$\nabla f \left(x_0,t\right)$ & Gradient of $f=f(x,t)$, with respect to $x$ variable,\\
 & for any $f $ such that $f(.,t)$ is differentiable at $x_0$\\
$\mathbf{\dot{X}}$ & Derivative with respect to $t$ of function $\mathbf{X}$, \\
& for any $ \mathbf{X}: [0,T] \rightarrow \mathbb{R}^{p}$ \\
$\partial_{i}$ & Partial derivative operator with respect to $x_i$ variable\\
$\partial_{t}$ & Partial derivative operator with respect to $t$ variable\\
$\langle , \rangle$ & Dot product  \\
$C^{k}\left( \Omega; \mathbb{R}^{p}\right)$ &  Class of $C^k$ functions $f: \Omega \rightarrow \mathbb{R}^{p}$, $\Omega \subset \mathbb{R}^{N}$ open\\ & $k=0$ is the class of continuos functions. We denote $C^0=C$\\
$C^{k,1}\left( \Omega \times [0,T]; \mathbb{R}^{p} \right)$ &  Class of functions $f: \Omega \times [0,T] \rightarrow \mathbb{R}^{p}$, $\Omega \subset \mathbb{R}^{N}$ open,\\
& $C^k$ in $\Omega$ and $C^1$ in [0,T] \\
$(f \vee g) (x)$ & Maximum value between $f(x)$ and $g(x)$ \\
\end{tabular}
\end{table}

\end{document}